\documentclass[12pt,fullpage]{article}

\usepackage{amsfonts}
\usepackage{xr}
\usepackage{graphicx}
\usepackage{amsmath}
\usepackage{epsfig}
\input epsf

\def\Z{\mathbb{Z}}
\def\R{\mathbb{R}}
\def\N{\mathbb{N}}
\def\tilde{\widetilde}
\def\epsilon{\varepsilon}
\def\trait (#1) (#2) (#3){\vrule width #1pt height #2pt depth #3pt}
\def\fin{\hfill\trait (0.1) (5) (0) \trait (5) (0.1) (0) \kern-5pt \trait (5) (5) (-4.9) \trait (0.1) (5) (0)}

\newcommand{\SE}{\setcounter{equation}{0} \section}
\newcommand{\be}{\begin{equation}}
\newcommand{\ee}{\end{equation}}
\newcommand{\baa}{\begin{array}}
\newcommand{\eaa}{\end{array}}
\newcommand{\ba}{\begin{eqnarray}}
\newcommand{\ea}{\end{eqnarray}}

\newtheorem{theo}{\bf Theorem}[section]
\newtheorem{lem}[theo]{\bf Lemma}
\newtheorem{pro}[theo]{\bf Proposition}
\newtheorem{cor}[theo]{\bf Corollary}
\newtheorem{defi}[theo]{\bf Definition}
\newtheorem{rem}[theo]{\bf Remark}

\begin{document}
\date{}
\title{\bf{Generalized transition waves and their properties}}
\author{Henri Berestycki$^{\hbox{\small{ a}}}$ and Fran\c cois Hamel$^{\hbox{\small{ b}}}$\thanks{Both authors are supported by the French ``Agence Nationale de la Recherche" within the project PREFERED. The second author is also indebted to the Alexander von~Humboldt Foundation for its support.}\\
\\
\footnotesize{$^{\hbox{a }}$EHESS, CAMS, 54 Boulevard Raspail, F-75006 Paris, France}\\
\footnotesize{\& University of Chicago, Department of Mathematics}\\
\footnotesize{5734 S. University Avenue, Chicago, IL 60637, USA}\\
\footnotesize{$^{\hbox{b }}$Aix-Marseille Universit\'e and Institut Universitaire de France}\\
\footnotesize{LATP, Avenue Escadrille Normandie-Niemen, F-13397 Marseille Cedex 20, France}}
\maketitle

\begin{abstract} In this paper, we generalize the usual notions of waves, fronts and propagation speeds in a very general setting. These new notions, which cover all usual situations, involve uniform limits, with respect to the geodesic distance, to a family of hypersurfaces which are parametrized by time. We prove the existence of new such waves for some time-dependent reaction-diffusion equations, as well as general intrinsic properties, some monotonicity properties and some uniqueness results for almost planar fronts. The classification results, which are obtained under some appropriate assumptions, show the robustness of our general definitions.
\end{abstract}

\tableofcontents


\SE{Introduction and main results}\label{intro}

We first introduce the definition of generalized transition waves and state some intrinsic general properties. We then give some specifications, including the notion of global mean speed, as well as some standard and new examples of such waves. Lastly, we state some of their important qualitative pro\-perties. We complete this section with some further possible extensions.


\subsection{Definition of generalized transition waves}

Traveling fronts describing the transition between two different states are a special important class of time-global solutions of evolution partial differential equations. One of the simplest examples is concerned with the homogeneous scalar semilinear parabolic equation
\begin{equation}\label{homo}
u_t=\Delta u+f(u)\hbox{ in }{\mathbb R}^N,
\end{equation}
where $u=u(t,x)$ and $\Delta$ is the Laplace operator with respect to the spatial variables in $\R^N$. In this case, assuming $f(0)=f(1)=0$, a planar traveling front connecting the uniform steady states $0$ and $1$ is a solution of the type
$$u(t,x)=\phi(x\cdot e-ct)$$
such that $\phi:\R\to[0,1]$ satisfies $\phi(-\infty)=1$ and $\phi(+\infty)=0$. Such a solution propagates in a given unit direction $e$ with the speed $c$. Existence and possible uniqueness of such fronts, formul{\ae} for the speed(s) of propagation are well-known \cite{aw,f,kpp} and depend upon the profile of the function $f$ on $[0,1]$.\par
In this paper, we generalize the standard notion of traveling fronts. That will allow us to consider new situations, that is new geometries or more complex equations. We provide explicit examples of new types of waves and we prove some qualitative properties. Although the definitions given below hold for general evolution equations (see Section~\ref{sec14}), we mainly focus on parabolic problems, that is we consider reaction-diffusion-advection equations, or systems of equations, of the type
\begin{equation}\label{eq}\left\{\begin{array}{l}
u_t=\nabla_x\cdot(A(t,x)\nabla_x u)+q(t,x)\cdot\nabla_x u+f(t,x,u)\hbox{ in }\Omega,\vspace{5pt}\\
g[t,x,u]=0\hbox{ on }\partial\Omega,\end{array}\right.
\end{equation}
where the unknown function $u$, defined in ${\mathbb R}\times\overline{\Omega}$, is in general a vector field
$$u=(u_1,\cdots,u_m)\in{\mathbb R}^m$$
and $\Omega$ is a globally smooth non-empty open connected subset of ${\mathbb R}^N$ with outward unit normal vector field $\nu$. By globally smooth, we mean that there exists $\beta>0$ such that $\Omega$ is globally of class $C^{2,\beta}$, that is there exist $r_0>0$ and $M>0$ such that, for all $y\in\partial\Omega$, there is a rotation $R_y$ of $\R^N$ and there is a $C^{2,\beta}$ map $\phi_y:\overline{B}^{N-1}_{2r_0}\to\R$ such that $\phi_y(0)=0$, $\|\phi_y\|_{C^{2,\beta}\big(\overline{B}^{N-1}_{2r_0}\big)}\le M$ and
$$\Omega\cap B(y,r_0)\!=\!\left[y+R_y\!\left(\{x\in\R^N;\ \!\!(x_1,\ldots,x_{N-1})\!\in\! \overline{B}^{N-1}_{2r_0},\ \!\!\phi_y(x_1,\ldots,x_{N-1})\!<\!x_N\}\right)\right]\!\cap\! B(y,r_0),$$
where $B(y,r_0)=\{x\in\R^N;\ |x-y|<r_0\}$, $|\ |$ denotes the Euclidean norm in $\R^N$ and, for any $s>0$, $\overline{B}^{N-1}_s$ is the closed Euclidean ball of $\R^{N-1}$ with center~$0$ and radius $s$ (notice in particular that~$\R^N$ is globally smooth).\par
Let us now list the general assumptions on the coefficients of (\ref{eq}). The diffusion matrix field
$$(t,x)\mapsto A(t,x)=(a_{ij}(t,x))_{1\le i,j\le N}$$
is assumed to be of class $C^{1,\beta}({\mathbb R}\times\overline{\Omega})$ and there exist $0<\alpha_1\le\alpha_2$ such that
$$\alpha_1|\xi|^2\le a_{ij}(t,x)\xi_i\xi_j\le\alpha_2|\xi|^2\hbox{ for all }(t,x)\in{\mathbb R}\times\overline{\Omega}\hbox{ and }\xi\in{\mathbb R}^N,$$
under the usual summation convention of repeated indices. The vector field
$$(t,x)\mapsto q(t,x)$$
ranges in ${\mathbb R}^N$ and is of class $C^{0,\beta}({\mathbb R}\times\overline{\Omega})$. The function $f:\R\times\overline{\Omega}\times\R^m\to\R^m$
$$(t,x,u)\mapsto f(t,x,u)$$
is assumed to be of class $C^{0,\beta}$ in $(t,x)$ locally in $u\in{\mathbb R}^m$, and locally Lipschitz-continuous in~$u$, uniformly with respect to $(t,x)\in{\mathbb R}\times\overline{\Omega}$. Lastly, the boundary conditions
$$g[t,x,u]=0\hbox{ on }\partial\Omega$$
may for instance be of the Dirichlet, Neumann, Robin or tangential types, or may be nonlinear or heterogeneous as well. The notation $g[t,x,u]=0$ means that this condition may involve not only $u(t,x)$ itself but also other quantities depending on $u$, like its derivatives for instance.\par
Throughout the paper, $d_{\Omega}$ denotes the geodesic distance in $\overline{\Omega}$, that is, for every pair~$(x,y)\in\overline{\Omega}\times\overline{\Omega}$, $d_{\Omega}(x,y)$ is the infimum of the arc lengths of all $C^1$ curves joining $x$ to $y$ in $\overline{\Omega}$. We assume that~$\Omega$ has an infinite diameter with respect to the geodesic distance~$d_{\Omega}$, that is $\hbox{diam}_{\Omega}(\Omega)=+\infty$, where
$$\hbox{diam}_{\Omega}(E)=\sup\big\{d_{\Omega}(x,y);\ (x,y)\in E\times E\big\}$$
for any $E\subset\Omega$. For any two subsets $A$ and $B$ of $\overline{\Omega}$, we set
$$d_{\Omega}(A,B)=\inf\big\{d_{\Omega}(x,y);\ (x,y)\in A\times B\big\}.$$
For $x\in\overline{\Omega}$ and $r>0$, we set
$$B_{\Omega}(x,r)=\big\{y\in\overline{\Omega};\ d_{\Omega}(x,y)<r\big\}\ \hbox{ and }\ S_{\Omega}(x,r)=\big\{y\in\overline{\Omega};\ d_{\Omega}(x,y)=r\big\}.$$\par
The following definition of a generalized transition wave, which has a geometric essence, involves two families $(\Omega^-_t)_{t\in{\mathbb R}}$ and $(\Omega^+_t)_{t\in{\mathbb R}}$ of open nonempty and unbounded subsets of $\Omega$ such that
\begin{equation}\label{omegapm}\left\{\baa{l}
\Omega^-_t\cap\Omega^+_t=\emptyset,\vspace{5pt}\\
\partial\Omega^-_t\cap\Omega=\partial\Omega^+_t\cap\Omega=:\Gamma_t,\vspace{5pt}\\
\Omega^-_t\cup\Gamma_t\cup\Omega^+_t=\Omega,\vspace{5pt}\\
\sup\big\{d_{\Omega}(x,\Gamma_t);\ x\in\Omega^+_t\big\}=+\infty,\vspace{5pt}\\
\sup\big\{d_{\Omega}(x,\Gamma_t);\ x\in\Omega^-_t\big\}=+\infty,\eaa\right.
\end{equation}
for all $t\in\R$. In other words, $\Gamma_t$ splits $\Omega$ into two parts, namely $\Omega^-_t$ and $\Omega^+_t$ (see Figure~\ref{fig1} below). The unboundedness of the sets $\Omega^{\pm}_t$ means that these sets have infinite diameters with respect to geodesic distance~$d_{\Omega}$. Moreover, for each~$t\in\R$, these sets are assumed to contain points which are as far as wanted from the interface~$\Gamma_t$. We further impose that
\be\label{unifgamma}
\sup\Big\{d_{\Omega}(y,\Gamma_t);\ y\in\overline{\Omega^{\pm}_t}\ \cap\  S_{\Omega}(x,r)\Big\}\to+\infty\hbox{ as }r\to+\infty\hbox{ uniformly in }t\in\R,\ x\in\Gamma_t
\ee
and that the interfaces $\Gamma_t$ are made of a finite number of graphs. By the latter we mean that, when $N\ge 2$, there is an integer $n\ge 1$ such that, for each $t\in\R$, there are $n$ open subsets~$\omega_{i,t}\subset\R^{N-1}$, $n$ continuous maps $\psi_{i,t}:\omega_{i,t}\to\R$ and $n$ rotations $R_{i,t}$ of~$\R^N$ $($for all~$1\le i\le n)$, such that
\be\label{omegapmbis}
\Gamma_t\subset\mathop{\bigcup}_{1\le i\le n}R_{i,t}\big(\big\{x\in\R^N;\ (x_1,\ldots,x_{N-1})\in\omega_{i,t},\ x_N=\psi_{i,t}(x_1,\ldots,x_{N-1})\big\}\big).
\ee
In dimension $N=1$, the above condition reduces to the existence of an integer $n\ge 1$ such that $\Gamma_t$ is made of at most $n$ points, that is $\Gamma_t=\{x^1_t,\ldots,x^n_t\}$ for each $t\in\R$ (where the real numbers $x^i_t$ may not be all pairwise distinct). As far as the condition~(\ref{unifgamma}) is concerned, its exact definition is
$$\left\{\baa{l}
\inf\Big\{\sup\big\{d_{\Omega}(y,\Gamma_t);\ y\in\overline{\Omega^+_t}\cap S_{\Omega}(x,r)\big\};\ t\in\R,\ x\in\Gamma_t\Big\}\to+\infty\vspace{5pt}\\
\inf\Big\{\sup\big\{d_{\Omega}(y,\Gamma_t);\ y\in\overline{\Omega^-_t}\cap S_{\Omega}(x,r)\big\};\ t\in\R,\ x\in\Gamma_t\Big\}\to+\infty\eaa\right.\hbox{ as }r\to+\infty.$$
It means that, for every point~$x\in\Gamma_t$, there are some points in both $\Omega^+_t$ and $\Omega^-_t$ which are far from $\Gamma_t$ and are at the same distance~$r$ from~$x$, when $r$ is large. The reason why this condition is used will become clearer in the following definition of transition waves.

\begin{defi}\label{def1} {\rm{(Generalized transition wave)}} Let $p^{\pm}:\R\times\overline{\Omega}\to\R^m$ be two classical solutions of $(\ref{eq})$. A {\rm{(generalized) transition wave}} connecting $p^-$ and $p^+$ is a time-global classical\footnote{Actually, from standard parabolic interior estimates, any classical solution of (\ref{eq}) is such that $u$, $u_t$,~$u_{x_i}$ and $u_{x_ix_j}$, for all $1\le i,j\le N$, are locally H\"older continuous in $\R\times\Omega$.} solution $u$ of $(\ref{eq})$ such that $u\not\equiv p^{\pm}$ and there exist some sets $(\Omega^{\pm}_t)_{t\in\R}$ and $(\Gamma_t)_{t\in\R}$ satisfying~$(\ref{omegapm})$,~$(\ref{unifgamma})$ and~$(\ref{omegapmbis})$ with
\be\label{defunif}
u(t,x)-p^{\pm}(t,x)\to 0\hbox{ uniformly in }t\in\R\hbox{ as }d_{\Omega}(x,\Gamma_t)\to+\infty\hbox{ and }x\in\overline{\Omega^{\pm}_t},
\ee
that is, for all $\epsilon>0$, there exists $M$ such that
$$\forall\, t\in\R,\ \forall\,x\in\overline{\Omega^{\pm}_t},\quad\big( d_{\Omega}(x,\Gamma_t)\ge M\big)\Longrightarrow\big(|u(t,x)-p^{\pm}(t,x)|\le\epsilon\big).$$
\end{defi}

Let us comment with words the key point in the above Definition~\ref{def1}.\footnote{Definition~\ref{def1} of generalized transition waves is slightly more precise than the one used in our companion paper~\cite{bh4}. In the present paper, we impose in the definition itself additional geometric conditions on the sets~$(\Omega^{\pm}_t)_{t\in\R}$ and~$(\Gamma_t)_{t\in\R}$, the meaning of which is explained in this paragraph.} Namely, a central role is played by the uniformity of the limits
$$u(t,x)-p^{\pm}(t,x)\to 0$$
as $d_{\Omega}(x,\Gamma_t)\to+\infty$ and $x\in\overline{\Omega^{\pm}_t}$. These limits hold far away from the hypersurfaces $\Gamma_t$ inside~$\Omega$. To make the definition meaningful, the distance which is used is the distance geodesic~$d_{\Omega}$. It is the right notion to fit with the geometry of the underlying domain. Furthermore, it is necessary to describe the propagation of transition waves in domains such as curved cylinders (like in the joint figure), spiral-shaped domains, exterior domains, etc. Roughly speaking, these limiting conditions~(\ref{defunif}), together with~(\ref{unifgamma}) and~(\ref{omegapmbis}), mean that the transition between the limiting states~$p^-$ and~$p^+$ is made of a finite number of neighborhoods of graphical interfaces, the width of these neighborhoods being bounded uniformly in time. Therefore, the region where a transition wave~$u$ connecting~$p^-$ and~$p^+$ is not close to~$p^{\pm}$ has a uniformly bounded width. This is the reason why the word ``transition", referring to the intuitive notion of spatial transition, is used to give a name to the objects introduced in Definition~\ref{def1}.\par
\begin{figure}\label{fig1}\centering
\includegraphics[width=7cm]{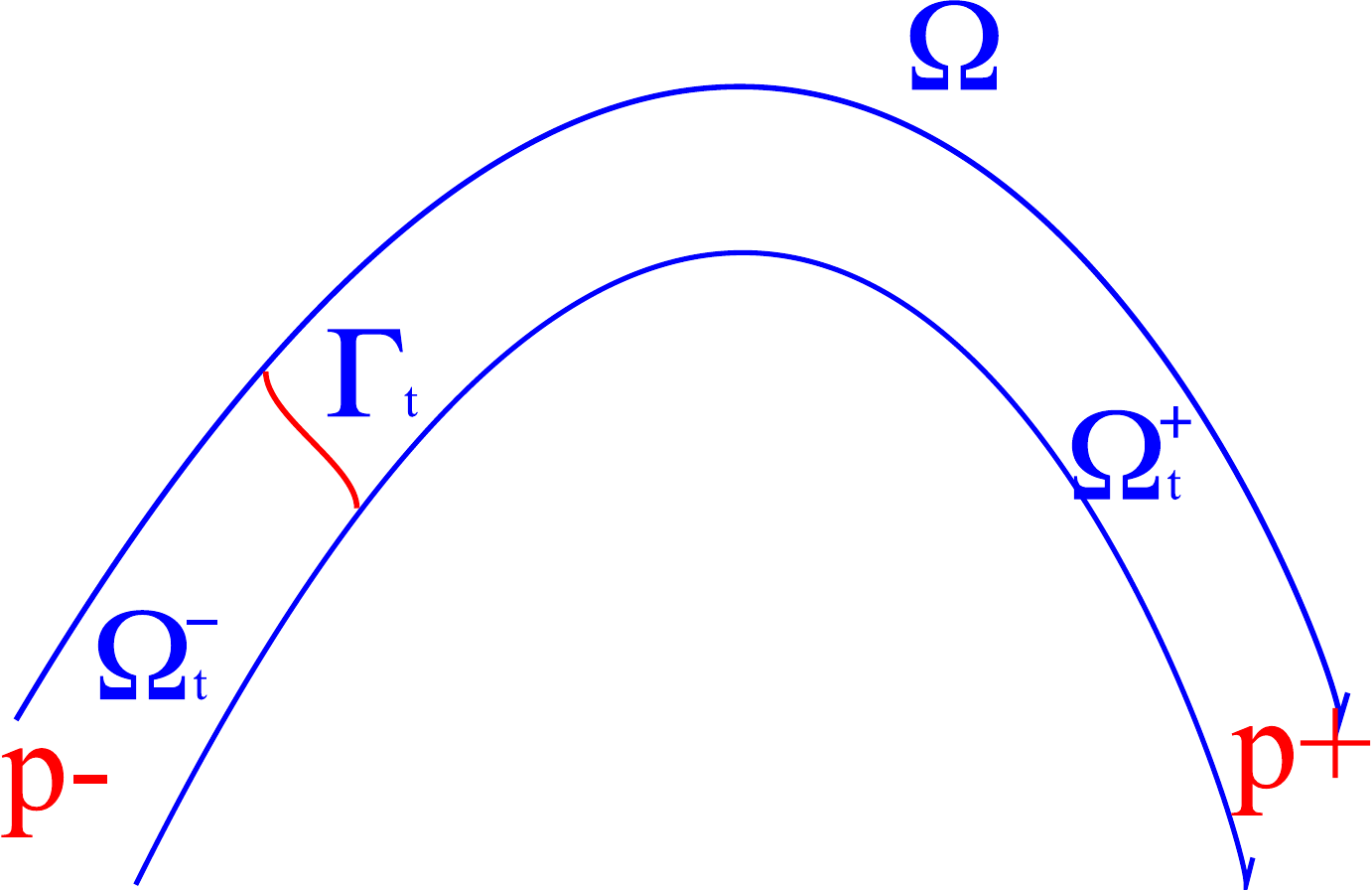}
\caption{A schematic picture of the sets $\Omega^{\pm}_t$ and $\Gamma_t$}
\end{figure}
We point out that, in Definition~\ref{def1}, the limiting states $p^{\pm}$ of a transition wave~$u$ are imposed to solve~(\ref{eq}). In other words, a transition wave is by definition a spatial connection between two other solutions. Thus, if $\epsilon^{\pm}$ are any two functions defined in $\R\times\overline{\Omega}$ such that~$\epsilon^{\pm}(t,x)\to 0$ as $d_{\Omega}(x,\Gamma_t)\to+\infty$ and $x\in\overline{\Omega^{\pm}_t}$ uniformly in $t$, and if $u$ is any time-global solution of~(\ref{eq}), then $u$ is in general not a transition wave between $u+\epsilon^-$ and $u+\epsilon^+$, because the limiting states $u+\epsilon^{\pm}$ do not solve~(\ref{eq}) in general. The requirement that the limiting states~$p^{\pm}$ of a transition wave~$u$ solve~(\ref{eq}) is then made in order to avoid the introduction of artificial and useless objects.\par
In Definition~\ref{def1}, the sets $(\Omega^{\pm}_t)_{t\in\R}$ and $(\Gamma_t)_{t\in\R}$ are not uniquely determined, given a generalized transition wave. Nevertheless, in the scalar case, under some assumptions on $p^{\pm}$ and $\Omega^{\pm}_t$ and oblique Neumann boundary conditions on $\partial\Omega$, the sets $\Gamma_t$ somehow reflect the location of the level sets of $u$. Namely, the following result holds:

\begin{theo}\label{pro1} Assume that $m=1$ $($scalar case$)$, that $p^{\pm}$ are constant solutions of $(\ref{eq})$ such that $p^-<p^+$ and let $u$ be a time-global classical solution of $(\ref{eq})$ such that
$$\big\{u(t,x);\ (t,x)\in\R\times\overline{\Omega}\big\}=(p^-,p^+)$$
and
$$g[t,x,u]=\mu(t,x)\cdot\nabla_xu(t,x)=0\ \hbox{ on }\R\times\partial\Omega,$$
for some unit vector field $\mu\in C^{0,\beta}(\R\times\partial\Omega)$ such that
$$\inf\big\{\mu(t,x)\cdot\nu(x);\ (t,x)\in\R\times\partial\Omega\big\}>0.\footnote{Therefore, $u$ and its derivatives $u_t$, $u_{x_i}$ and $u_{x_ix_j}$, for all $1\le i,j\le N$, are bounded and globally H\"older continuous in $\R\times\overline{\Omega}$.}$$\par
{\rm{1.}} Assume that $u$ is a generalized transition wave connecting $p^-$ and $p^+$, or $p^+$ and $p^-$, in the sense of Definition~$\ref{def1}$ and that there exists $\tau>0$ such that
\be\label{tau}
\sup\big\{d_{\Omega}(x,\Gamma_{t-\tau});\ t\in\R,\ x\in\Gamma_t\big\}<+\infty.
\ee
Then
\be\label{lambda}
\forall\,\lambda\in(p^-,p^+),\ \ \sup\big\{d_{\Omega}(x,\Gamma_t);\ u(t,x)=\lambda\big\}<+\infty
\ee
and
\be\label{infsup}
\forall\,C\ge 0,\ \ p^-<\inf\big\{u(t,x);\ d_{\Omega}(x,\Gamma_t)\!\le\!C\big\}\le\sup\big\{u(t,x);\ d_{\Omega}(x,\Gamma_t)\!\le\!C\big\}<p^+.
\ee\par
{\rm{2.}} Conversely, if~$(\ref{lambda})$ and~$(\ref{infsup})$ hold for some choices of sets $(\Omega^{\pm}_t,\Gamma_t)_{t\in\R}$ satis\-fying~$(\ref{omegapm})$,~$(\ref{unifgamma})$ and~$(\ref{omegapmbis})$, and if there is $d_0>0$ such that the sets
$$\big\{(t,x)\in{\mathbb{R}}\times\overline{\Omega};\ x\in\overline{\Omega^+_t},\ d_{\Omega}(x,\Gamma_t)\ge d\big\}$$
and
$$\big\{(t,x)\in{\mathbb{R}}\times\overline{\Omega};\ x\in\overline{\Omega^-_t},\ d_{\Omega}(x,\Gamma_t)\ge d\big\}$$
are connected for all $d\ge d_0$, then $u$ is a generalized transition wave connecting $p^-$ and $p^+$, or $p^+$ and $p^-$.
\end{theo}

The assumption (\ref{tau}) means that the interfaces $\Gamma_t$ and $\Gamma_{t-\tau}$ are in some sense not too far from each other. For instance, if all $\Gamma_t$ are parallel hyperplanes in $\Omega=\R^N$, then the assumption (\ref{tau}) means that the distance between $\Gamma_t$ and $\Gamma_{t-\tau}$ is bounded independently of~$t$, for some $\tau>0$. As far as the connectedness assumptions made in part~2 of Theorem~\ref{pro1} are concerned, they are a topological ingredient in the proof, to guarantee the uniform convergence of $u$ to $p^{\pm}$ or $p^{\mp}$ far away from $\Gamma_t$ in $\Omega^{\pm}_t$.


\subsection{Some specifications and the notion of global mean speed}

In this section, we define the more specific notions of fronts, pulses, invasions (or traveling waves) and almost planar waves, as well as the concept of global mean speed, when it exists. These notions are related to some analytical or geometric properties of the limiting states~$p^{\pm}$ or of the sets $(\Omega^{\pm}_t)_{t\in\R}$ and $(\Gamma_t)_{t\in\R}$, and are listed in the following definitions, where~$u$ denotes a transition wave connecting $p^-$ and $p^+$, associated to two families $(\Omega^{\pm}_t)_{t\in\R}$ and~$(\Gamma_t)_{t\in\R}$, in the sense of Definition~\ref{def1}.

\begin{defi}\label{frontpulse} {\rm{(Fronts and spatially extended pulses)}} Let $p^{\pm}=(p^{\pm}_1,\cdots,p^{\pm}_m)$. We say that the transition wave $u$ is a {\rm{front}} if, for each $1\le k\le m$, either
$$\inf\big\{p^+_k(t,x)-p^-_k(t,x);\ x\in\overline{\Omega}\big\}>0\hbox{ for all }t\in{\mathbb R}$$
or
$$\inf\big\{p^-_k(t,x)-p^+_k(t,x);\ x\in\overline{\Omega}\big\}>0\hbox{ for all }t\in{\mathbb R}.$$
The transition wave $u$ is a {\rm{spatially extended pulse}} if $p^{\pm}$ depend only on~$t$ and $p^-(t)=p^+(t)$ for all $t\in\R$.
\end{defi}

In the scalar case ($m=1$), our definition of a front corresponds to the natural extension of the usual notion of a front connecting two different constants. In the pure vector case ($m\ge 2$), if a bounded $C^{0,\beta}(\R\times\overline{\Omega})$ transition wave $u=(u_1,\ldots,u_m)$ is a front for problem
$$u_t=\nabla_x\cdot(A(t,x)\nabla_x u)+q(t,x)\cdot\nabla_x u+f(t,x,u)$$
in the sense of Definitions~\ref{def1} and~\ref{frontpulse}, if $u_k\not\equiv p^{\pm}_k$ for some $1\le k\le m$, then the function $u_k$ is a front connecting $p^-_k$ and $p^+_k$ for the problem
$$(u_k)_t=\nabla_x\cdot(A(t,x)\nabla_xu_k)+q(t,x)\cdot\nabla_xu_k+\tilde{f}_k(t,x,u_k)$$
associated with the same sets $(\Omega^{\pm}_t)_{t\in\R}$ and $(\Gamma_t)_{t\in\R}$ as $u$, where
$$\tilde{f}_k(t,x,s)=f(t,x,u_1(t,x),\ldots,u_{k-1}(t,x),s,u_{k+1}(t,x),\ldots,u_m(t,x))$$
and $f=(f_1,\ldots,f_m)$. The same observation is valid for spatially extended pulses as well.

\begin{defi}\label{invasion} {\rm{(Invasions)}} We say that $p^+$ {\rm{invades}} $p^-$, or that $u$ is an {\rm{invasion}} of $p^-$ by $p^+$ $($resp. $p^-$ invades $p^+$, or $u$ is an invasion of $p^+$ by $p^-)$ if
$$\Omega^+_t\supset\Omega^+_s\ (\hbox{resp. }\Omega^-_t\supset\Omega^-_s)\hbox{ for all }t\ge s$$
and
$$d_{\Omega}(\Gamma_t,\Gamma_s)\to+\infty\hbox{ as }|t-s|\to+\infty.$$
\end{defi}

Therefore, if $p^+$ invades $p^-$ (resp. $p^-$ invades $p^+$), then $u(t,x)-p^{\pm}(t,x)\to 0$ as $t\to\pm\infty$ (resp. as $t\to\mp\infty$) locally uniformly in $\overline{\Omega}$ with respect to the distance $d_{\Omega}$. One can then say that, roughly speaking, invasions correspond to the usual idea of traveling waves. Notice that a generalized transition wave can always be viewed as a spatial connection between~$p^-$ and~$p^+$, while an invasion wave can also be viewed as a temporal connection between the limiting states $p^-$ and $p^+$.

\begin{defi} {\rm{(Almost planar waves in the direction $e$)}} We say that the generalized transition wave $u$ is {\rm{almost planar}} $($in the direction $e\in{\mathbb{S}}^{N-1})$ if, for all $t\in\mathbb{R}$, the sets $\Omega^{\pm}_t$ can be chosen so that
$$\Gamma_t=\big\{x\in\Omega;\ x\cdot e=\xi_t\big\}$$
for some $\xi_t\in\mathbb{R}$.
\end{defi}

By extension, we say that the generalized transition wave $u$ is {\it{almost planar}} in a moving direction $e(t)\in{\mathbb{S}}^{N-1}$ if, for all $t\in\mathbb{R}$, $\Omega^{\pm}_t$ can be chosen so that
$$\Gamma_t=\big\{x\in\Omega;\ x\cdot e(t)=\xi_t\big\}$$
for some $\xi_t\in\mathbb{R}$.\par
As in the usual cases (see Section~\ref{sec12}), an important notion which is attached to a generalized transition wave is that of its global mean speed of propagation, if any.

\begin{defi}\label{defspeed} {\rm{(Global mean speed of propagation)}} We say that a generalized transition wave $u$ associated to the families $(\Omega^{\pm}_t)_{t\in\R}$ and $(\Gamma_t)_{t\in\R}$ has {\rm{global mean speed}} $c\ (\ge 0)$ if
$$\frac{d_{\Omega}(\Gamma_t,\Gamma_s)}{|t-s|}\to c\ \hbox{ as }|t-s|\to+\infty.$$
We say that the transition wave $u$ is {\rm{almost-stationary}} if it has global mean speed $c=0$. We say that $u$ is {\rm{quasi-stationary}} if
$$\sup\ \{d_{\Omega}(\Gamma_t,\Gamma_s);\ (t,s)\in{\mathbb R}^2\}<+\infty,$$
and we say that $u$ is {\rm{stationary}} if it does not depend on $t$.
\end{defi}

The global mean speed $c$, if it exists, is unique. Moreover, under some reasonable assumptions, the global mean speed is an intrinsic notion, in the sense that it does not depend on the families $(\Omega^{\pm}_t)_{t\in\R}$ and $(\Gamma_t)_{t\in\R}$. This is indeed seen in the following result:

\begin{theo}\label{pro2}
In the general vectorial case $m\ge 1$, let $p^{\pm}$ be two solutions of $(\ref{eq})$ satisfying
$$\inf\big\{|p^-(t,x)-p^+(t,x)|;\ (t,x)\in{\mathbb R}\times\overline{\Omega}\big\}>0.$$
Let $u$ be a transition wave connecting $p^-$ and $p^+$ with a choice of sets $(\Omega^{\pm}_t)_{t\in\R}$ and $(\Gamma_t)_{t\in\R}$, satisfying $(\ref{omegapm})$,~$(\ref{unifgamma})$ and~$(\ref{omegapmbis})$. If $u$ has global mean speed $c$, then, for any other choice of sets $(\tilde{\Omega}^{\pm}_t)_{t\in\R}$ and $(\tilde{\Gamma}_t)_{t\in\R}$, satisfying $(\ref{omegapm})$,~$(\ref{unifgamma})$ and~$(\ref{omegapmbis})$,~$u$ has a global mean speed and this global mean speed is equal to~$c$.
\end{theo}


\subsection{Usual cases and new examples}\label{sec12}

In this subsection, we list some basic examples of transition waves, which correspond to the classical notions in the standard situations. We also state the existence of new examples of transition fronts in a time-dependent medium.\par
For the homogeneous equation (\ref{homo}) in ${\mathbb R}^N$, a solution
$$u(t,x)=\phi(x\cdot e-ct),$$
with $\phi(-\infty)=1$ and $\phi(+\infty)=0$ (assuming $f(0)=f(1)=0$) is an (almost) planar invasion front connecting $p^-=1$ and $p^+=0$, with (global mean) speed $|c|$. The uniform stationary state $p^-=1$ (resp. $p^+=0$) invades the uniform stationary $p^+=0$ (resp. $p^-=1$) if $c>0$ (resp. $c<0$). The sets $\Omega^{\pm}_t$ can for instance be defined as
$$\Omega^{\pm}_t=\big\{x\in{\mathbb R}^N;\ \pm(x\cdot e-ct)>0\big\}.$$
The general definitions that we just gave also generalize the classical notions of pulsating traveling fronts in spatially periodic media (see \cite{bh2,bhn1,bhn2,bhr2,h,hr2,hz,skt,x1,x2}) with possible periodicity or almost-periodicity in time (see \cite{gf,nad,nrx,nx,sh1,sh2,sh3}) or in spatially recurrent media (see \cite{m}).\par
We point out that the limiting states~$p^{\pm}(t,x)$ are not assumed to be constant in general. It is indeed important to let the possibility of transition waves connecting time- or space-dependent limiting states. In the aforementioned references in the periodic case, the limiting states are typically periodic as well. Let us mention here another situation, corresponding to a one-dimensional medium which is asymptotically homogeneous but not uniformly homogeneous, and let us explain what a transition wave can be in this case. Namely, consider an equation of the type
$$u_t=u_{xx}+f(x,u),$$
where $u:\R\times\R\to\R^m$, $f(x,a_1)=0$ for all $x\in\R$, $f(x,a_2)\to0$ as $x\to-\infty$, $f(x,a_3)\to 0$ as $x\to+\infty$, and $a_1$,~$a_2$ and~$a_3$ are three distinct vectors in~$\R^m$. The homogeneous states~$a_2$ and~$a_3$ are solutions of the limiting equations obtained as $x\to-\infty$ and $x\to+\infty$ respectively, but these states do not solve the original equation in general since $f(x,a_2)$ and $f(x,a_3)$ are not identically equal to $0$ in general. One can then wonder what could be a generalized transition wave~$u(t,x)$ connecting~$p^-=a_1$ to another limiting state~$p^+$, with a single interface~$\Gamma_t=\{x_t\}$ such that $\Omega^{\pm}_t=\{x\in\R;\ \pm(x-x_t)<0\}$ and $x_t\to\pm\infty$ as $t\to\pm\infty$. The limiting state~$p^+(t,x)$ such that $u(t,x)-p^+(t,x)\to 0$ as $x-x_t\to-\infty$ (uniformly in~$t\in\R$) cannot be~$a_2$ or~$a_3$ in general. A natural candidate could be a solution~$p^+(x)$ of the stationary equation
$${p^+}''(x)+f(x,p^+(x))=0,\ \ x\in\R,$$
such that $p^+(-\infty)=a_2$ and $p^+(+\infty)=a_3$. If such a solution~$p^+$ exists, a transition wave connecting~$p^-=a_1$ and~$p^+$ and satisfying $\lim_{t\to\pm\infty}x_t=\pm\infty$ would then be such that~$u(t,x)\to a_2$ (resp.~$u(t,x)\to a_1$) as $x\to-\infty$ (resp. $x\to+\infty$) locally uniformly in~$t\in\R$, but $u(t,x)\to p^+(x)\not\equiv a_2$ as~$t\to+\infty$ locally uniformly in $x\in\R$. Without going into further details here, this simple example already illustrates the wideness of Definition~\ref{def1} and the possibility of new objects connecting general non-constant limiting states.\par
In the particular one-dimensional case, when equation (\ref{eq}) is scalar and when the limi\-ting states $p^{\pm}$ are ordered, say $p^+>p^-$, Definition~\ref{def1} corresponds to that of ``wave-like" solutions given in \cite{sh4}. However, Definition~\ref{def1} also includes more general situations involving complex heterogeneous geometries or media. Existence, uniqueness and stability results of generalized almost planar transition fronts in one-dimensional media or straight higher-dimensional cylinders with combustion-type nonlinearities and arbitrary spatial dependence have just been proved in \cite{mnrr,mrs,nr,z}. In general higher-dimensional domains, generalized transition waves which are not almost planar can also be covered by Definition~\ref{def1}: such transition waves are known to exist for the homogeneous equation (\ref{homo}) in $\R^N$ for usual types of nonlinearities $f$ (combustion, bistable, Kolmogorov-Petrovsky-Piskunov type), see \cite{bh4,boh,hmr,hn,hs,hu,nt,t1,t2} for details. Further on, other situations can also be investigated, such as the case when some coefficients of (\ref{eq}) are locally perturbed and more complex geometries, like exterior domains (the existence of almost planar fronts in exterior domains with bistable nonlinearity $f$ has just been proved in \cite{bhm}), curved cylinders, spirals, etc can be considered.\par
It is worth mentioning that, even in dimension~1, Definition~\ref{def1} also includes a very interesting class of transition wave solutions which are known to exist and which do not fall within the usual notions, that is invasion fronts which have no specified global mean speed. For instance, for (\ref{homo}) in dimension $N=1$, if $f=f(u)$ satisfies
\be\label{fconcave}
f\hbox{ is }C^2\hbox{ concave in }[0,1],\hbox{ positive in }(0,1)\hbox{ and }f(0)=f(1)=0,
\ee
then there are invasion fronts connecting $0$ and $1$ for which $\Omega^-_t=(x_t,+\infty)$, $\Omega^+_t=(-\infty,x_t)$ and
$$\frac{x_t}{t}\to c_1\hbox{ as }t\to-\infty\hbox{ and }\frac{x_t}{t}\to c_2\hbox{ as }t\to+\infty$$ with $2\sqrt{f'(0)}\le c_1<c_2$ (see \cite{hn}). There are also some fronts for which $x_t/t\to c_1\ge 2\sqrt{f'(0)}$ as $t\to-\infty$ and $x_t/t\to+\infty$ as $t\to+\infty$. For further details, we refer to \cite{bh4,hn}.\par
In the companion survey paper~\cite{bh4}, we made a detailed presentation of the usual particular cases of transition waves covered by Definition~\ref{def1}. We explained and compared the notions of fronts which had been introduced earlier, starting from the simplest situations and going to the most general ones. In the present paper, in addition to the intrinsic properties of the generalized transition waves stated in Theorems~\ref{pro1} and~\ref{pro2} above, we mainly focus on the proof of some important qualitative properties, including some monotonicity and uniqueness results, and on the application of these qualitative properties in order to get Liouville-type results in some particular situations. In doing so, we prove that, under some assumptions, the generalized transition waves reduce to the standard traveling or pulsating fronts in homogeneous or periodic media. These qualitative properties are stated in the next subsection~\ref{secqual}. In a forthcoming paper, we deal with a general method to prove the existence of transition waves in a broad framework. However, in the present paper, in order to illustrate the interest of the above definitions, we also analyze a specific example which had not been considered in the literature. We prove the existence of new generalized transition waves, which in general do not have any global mean speed, for time-dependent equations. Namely, we consider one-dimensional reaction-diffusion equations of the type
\be\label{eqkpp}
u_t=u_{xx}+f(t,u)
\ee
where the function $f:\R\times[0,1]\to\R$ is of class $C^1$ and satisfies:
\be\label{ftime}\left\{\baa{l}
\forall\,t\in\R,\quad f(t,0)=f(t,1)=0,\vspace{5pt}\\
\forall\,(t,s)\in\R\times[0,1],\quad f(t,s)\ge 0,\vspace{5pt}\\
\exists\,t_1<t_2\in\R,\ \exists\,f_1,\,f_2\in C^1([0,1];\R),\vspace{3pt}\\
\qquad\qquad\forall\,(t,s)\in(-\infty,t_1]\times[0,1],\quad f(t,s)=f_1(s),\vspace{3pt}\\
\qquad\qquad\forall\,(t,s)\in[t_2,+\infty)\times[0,1],\quad f(t,s)=f_2(s),\vspace{3pt}\\
\qquad\qquad f'_1(0)>0,\ f'_2(0)>0,\vspace{3pt}\\
\qquad\qquad\forall\,s\in(0,1),\quad f_1(s)>0,\ f_2(s)>0.\eaa\right.
\ee
In other words, the function $f$ is time-independent and non-degenerate at $0$ for times less than~$t_1$ and larger than $t_2$, and for the times $t\in(t_1,t_2)$, the functions $f(t,\cdot)$ are just assumed to be nonnegative, but they may a priori vanish. If $f_1$ and $f_2$ are equal, then the nonlinearity~$f(t,s)$ can be viewed as a time-local perturbation of a time-independent equation. But, it is worth noticing that the functions $f_1$ and $f_2$ are not assumed to be equal nor even compared in general. When $t\le t_1$, classical traveling fronts
$$(t,x)\in\R^2\mapsto\varphi_{1,c}(x-ct)\in[0,1]$$
such that $\varphi_{1,c}(-\infty)=1$ and $\varphi_{1,c}(+\infty)=0$ are known to exist, for all and only all speeds $c\ge c^*_1$, where $c_1^*\ge 2\sqrt{f'_1(0)}>0$ only depends on $f_1$ (see e.g. \cite{aw}). The open questions are to know how these traveling fronts behave during the time interval $[t_1,t_2]$ and whether they can subsist and at which speed, if any, they travel after the time $t_2$. Indeed, it is also known that, when $t\ge t_2$, there exist classical traveling fronts
$$(t,x)\in\R^2\mapsto\varphi_{2,c}(x-ct)\in[0,1]$$
such that $\varphi_{2,c}(-\infty)=1$ and $\varphi_{2,c}(+\infty)=0$  for all and only all speeds $c\ge c^*_2$, where $c_2^*\ge 2\sqrt{f'_2(0)}>0$ only depends on $f_2$. The following result provides an answer to these questions and shows the existence of generalized transition waves connecting $0$ and $1$ for equation (\ref{eqkpp}), which fall within our general definitions and do not have any global mean speed in general. To state the result, we need a few notations. For each $c\ge c^*_1$, we set
\be\label{lambda1c}
\lambda_{1,c}=\left\{\baa{ll}
\displaystyle{\frac{c_1-\sqrt{c_1^2-4f'_1(0)}}{2}} & \hbox{if }c>c^*_1,\vspace{3pt}\\
\displaystyle{\frac{c^*_1+\sqrt{{c^*_1}^2-4f'_1(0)}}{2}} & \hbox{if }c=c^*_1.\eaa\right.
\ee
We also denote
$$\lambda_2^{*,-}=\displaystyle{\frac{c^*_2-\sqrt{{c^*_2}^2-4f'_2(0)}}{2}}.$$

\begin{theo}\label{thtime}
For equation $(\ref{eqkpp})$ under the assumption $(\ref{ftime})$, there exist transition invasion fronts connecting $p^-=0$ and $p^+=1$, for which $\Omega^{\pm}_t=\{x\in\R;\ \pm(x-x_t)<0\}$, $\Gamma_t=\{x_t\}$ for all $t\in\R$,
$$x_t=c_1t\hbox{ for }t\le t_1\hbox{ and }\frac{x_t}{t}\to c_2\hbox{ as }t\to+\infty,$$
where $c_1$ is any given speed in $[c^*_1,+\infty)$ and
\be\label{c2}
c_2=\left\{\baa{ll}
\lambda_{1,c_1}+\displaystyle{\frac{f'_2(0)}{\lambda_{1,c_1}}} & \hbox{if }\lambda_{1,c_1}<\lambda_2^{*,-},\vspace{3pt}\\
c^*_2 & \hbox{if }\lambda_{1,c_1}\ge\lambda_2^{*,-}.\eaa\right.
\ee
\end{theo}

When $f_1=f_2$, then $c^*_1=c^*_2$ and the transition fronts constructed in Theorem~\ref{thtime} are such that $c_1=c_2$, whence they have a global mean speed $c=c_1=c_2$ in the sense of Definition~\ref{defspeed}. When $f_1\le f_2$ (resp. $f_1\ge f_2$), then $c^*_1\le c^*_2$ (resp. $c^*_1\ge c^*_2$), the inequalities $c_1\le c_2$ (resp. $c_1\ge c_2$) always hold and, for $c_1$ large enough so that $\lambda_{1,c_1}<\lambda_2^{*,-}$, the inequalities $c_1< c_2$ (resp. $c_1>c_2$) are strict if $f'_1(0)\neq f'_2(0)$ (hence, these transition fronts do not have any global mean speed).\par
In the general case, acceleration and slow down may occur simultaneously, for  the same equation (\ref{eqkpp}) with the same function $f$, according to the starting speed $c_1$: for instance, there are examples of functions $f_1$ and $f_2$ for which
$$c_2>c_1\hbox{ for all }c_1>c^*_1,\hbox{ and }c_2<c_1\hbox{ for }c_1=c^*_1.$$
To do so, it is sufficient to choose $f_2$ of the Kolmogorov-Petrovsky-Piskunov type, that is $f_2(s)\le f'_2(0)s$ in $(0,1)$ whence $c^*_2=2\sqrt{f'_2(0)}=2\lambda_2^{*,-}$, and to choose $f_1$ in such a way that $f'_1(0)<f'_2(0)$ and $c^*_1>c^*_2$ (for instance, if $f_2$ is chosen as above, if $M>0$ is such that $\sqrt{2M}>c^*_2$ and if
$$f_1(s)\ge\frac{M}{\epsilon}\times(1-|x-1+\epsilon|)\ \hbox{ on }[1-2\epsilon,1]$$
for $\epsilon>0$ small enough, then $c^*_1>c^*_2$ for $\epsilon$ small enough, see \cite{bns}).\par
Lastly, it is worth noticing that, in Theorem~\ref{thtime}, the speed $c_2$ of the position $x_t$ at large time is determined only from $c_1$, $f_1$ and $f_2$, whatever the profile of $f$ between times $t_1$ and~$t_2$ may be.

\begin{rem}{\rm The solutions $u$ constructed in Theorem~\ref{thtime} are by definition {\it spatial} transition fronts connecting $1$ and $0$. Furthermore, it follows from the proof given in Section~\ref{sec3} that these transition fronts can also be viewed as {\it temporal} connections between a classical traveling front with speed $c_1$ for the nonlinearity $f_1$ and another classical traveling front, with speed $c_2$, for the nonlinearity $f_2$.}
\end{rem}


\subsection{Qualitative properties}\label{secqual}

We now proceed to some further qualitative properties of generalized transition waves. {\it Throughout this subsection, $m=1$, i.e. we work in the scalar case, and $u$ denotes transition wave connec\-ting~$p^-$ and $p^+$, for equation $(\ref{eq})$, associated with families $(\Omega^{\pm}_t)_{t\in\R}$ and~$(\Gamma_t)_{t\in\R}$ satisfying properties $(\ref{omegapm})$, $(\ref{unifgamma})$, $(\ref{omegapmbis})$ and $(\ref{tau})$.} We also assume that $u$ and $p^{\pm}$ are globally bounded in~$\R\times\overline{\Omega}$ and that
\be\label{mu}
\mu(x)\cdot\nabla_xu(t,x)=\mu(x)\cdot\nabla_xp^{\pm}(t,x)=0\ \hbox{ on }\R\times\partial\Omega,
\ee
where $\mu$ is a $C^{0,\beta}(\partial\Omega)$ unit vector field such that
$$\inf\big\{\mu(x)\cdot\nu(x);\ x\in\partial\Omega\big\}>0.$$\par
First, we establish a general property of monotonicity with respect to time. 

\begin{theo}\label{th1}
Assume that $A$ and $q$ do not depend on $t$, that $f$ and $p^{\pm}$ are nondecreasing in $t$ and that there is $\delta>0$ such that
\be\label{fdecreasing}
s\mapsto f(t,x,s)\hbox{ is nonincreasing in }(-\infty,p^-(t,x)+\delta]\hbox{ and }[p^+(t,x)-\delta,+\infty)
\ee
for all $(t,x)\in\R\times\overline{\Omega}$. If $u$ is an invasion of $p^-$ by $p^+$ with
\be\label{kappa}
\kappa:=\inf\big\{p^+(t,x)-p^-(t,x);\ (t,x)\in\mathbb{R}\times\overline{\Omega}\big\}>0,
\ee
then $u$ satisfies
\be\label{ineqs}
\forall\ (t,x)\in\R\times\overline{\Omega},\quad p^-(t,x)<u(t,x)<p^+(t,x).
\ee
and $u$ is increasing in time $t$.
\end{theo}

Notice that if (\ref{ineqs}) holds a priori and if $f$ is assumed to be nonincreasing in $s$ for $s$ in $[p^-(t,x),p^-(t,x)+\delta]$ and $[p^+(t,x)-\delta,p^+(t,x)]$ only, instead of $(-\infty,p^-(t,x)+\delta]$ and $[p^+(t,x)-\delta,+\infty)$, then the conclusion of Theorem \ref{th1} (strict monotonicity of $u$ in $t$) holds. The simplest case is when $f=f(u)$ only depends on $u$ and $p^{\pm}$ are constants and both stable, that is $f'(p^{\pm})<0$.\par
The monotonicity result stated in Theorem \ref{th1} plays an important role in the following uniqueness and comparison properties for almost planar fronts:

\begin{theo}\label{th2} Under the same conditions as in Theorem~$\ref{th1}$, assume furthermore that $f$ and $p^{\pm}$ are independent of $t$, that $u$ is almost planar in some direction $e\in\mathbb{S}^{N-1}$ and has global mean speed $c\ge 0$, with the stronger property that
\be\label{Gamma}
\sup\big\{\big|\,d_{\Omega}(\Gamma_t,\Gamma_s)-c|t-s|\,\big|;\ (t,s)\in\R^2\big\}<+\infty,
\ee
where
$$\Gamma_t=\big\{x\in\Omega;\ x\cdot e-\xi_t=0\big\}\hbox{ and }\Omega^{\pm}_t=\big\{x\in\Omega;\ \pm(x\cdot e-\xi_t)<0\big\}.$$
Let $\tilde{u}$ be another globally bounded invasion front of $p^-$ by $p^+$ for equation $(\ref{eq})$ and $(\ref{mu})$, associated with
$$\tilde{\Gamma}_t=\big\{x\in\Omega;\ x\cdot e-\tilde{\xi}_t=0\big\}\hbox{ and }\tilde{\Omega}^{\pm}_t=\big\{x\in\Omega;\ \pm(x\cdot e-\tilde{\xi}_t)<0\big\}$$
and having global mean speed ${\tilde{c}}\ge 0$ such that
$$\sup\big\{\,\big|\,d_{\Omega}({\tilde{\Gamma}}_t,{\tilde{\Gamma}}_s)-{\tilde{c}}|t-s|\,\big|;\ (t,s)\in\R^2\big\}<+\infty.$$
Then $c={\tilde{c}}$ and there is $($the smallest$)$ $T\in\mathbb{R}$ such that
$$\tilde{u}(t+T,x)\ge u(t,x)\hbox{ for all }(t,x)\in\R\times\overline{\Omega}.$$ Furthermore, there exists a sequence $(t_n,x_n)_{n\in\N}$ in $\R\times\overline{\Omega}$ such that
$$(d_{\Omega}(x_n,\Gamma_{t_n}))_{n\in\N}\hbox{ is bounded and }\tilde{u}(t_n+T,x_n)-u(t_n,x_n)\to 0\hbox{ as }n\to+\infty.$$
Lastly, either $\tilde{u}(t+T,x)>u(t,x)$ for all $(t,x)\in\R\times\overline{\Omega}$ or $\tilde{u}(t+T,x)=u(t,x)$ for all $(t,x)\in\R\times\overline{\Omega}$.
\end{theo}

This result shows the uniqueness of the global mean speed among a certain class of almost planar invasion fronts. It also says that any two such fronts can be compared up to shifts. In particular cases listed below, uniqueness holds up to shifts. However, this uniqueness property may not hold in general.

\begin{rem}{\rm Notice that property~(\ref{Gamma}) and the fact that $u$ is an invasion imply that the speed~$c$ is necessarily (strictly) positive.}
\end{rem}

As a corollary of Theorem~\ref{th2}, we now state a result which is important in that it shows that, at least under appropriate conditions on $f$, our definition does not introduce new objects in some classical situations: it reduces to pulsating traveling fronts in periodic media and to usual traveling fronts when there is translation invariance in the direction of propagation.

\begin{theo}\label{cor1} Under the conditions of Theorem~$\ref{th2}$, assume that $\Omega$, $A$, $q$, $f$, $\mu$ and $p^{\pm}$ are periodic in $x$, in that there are positive real numbers $L_1,\ldots,L_N>0$ such that, for every vector $k=(k_1,\ldots,k_N)\in L_1\Z\times\cdots\times L_N\Z$,
$$\left\{\baa{l}
\Omega+k=\Omega,\vspace{5pt}\\
A(x+k)=A(x),\ \!q(x+k)=q(x),\ \!f(x+k,\cdot)=f(x,\cdot),\ \!p^{\pm}(x+k)=p^{\pm}(x)\hbox{ for all }x\in\overline{\Omega},\vspace{5pt}\\
\mu(x+k)=\mu(x)\hbox{ for all }x\in\partial\Omega.\eaa\right.$$\par
{\rm{(i)}} Then $u$ is a pulsating front, namely
\be\label{ptf}
u\left(t+\frac{\gamma\ k\cdot e}{c},x\right)=u(t,x-k)\hbox{ for all }(t,x)\in\mathbb{R}\times\overline{\Omega}\hbox{ and }k\in L_1\mathbb{Z}\times\cdots\times L_N\mathbb{Z},
\ee
where $\gamma=\gamma(e)\ge 1$ is given by
\be\label{gamma}
\gamma(e)=\mathop{\lim}_{(x,y)\in\overline{\Omega}\times\overline{\Omega},\ (x-y)\,\parallel\,e,\ |x-y|\to+\infty}\frac{d_{\Omega}(x,y)}{|x-y|}.
\ee
Furthermore, $u$ is unique up to shifts in $t$.\par
{\rm{(ii)}} Under the additional assumptions that $e$ is one of the axes of the frame, that $\Omega$ is invariant in the direction $e$ and that $A$, $q$, $f$, $\mu$ and $p^{\pm}$ are independent of $x\cdot e$, then $u$ actually is a classical traveling front, that is:
$$u(t,x)=\phi(x\cdot e-ct,x')$$
for some function $\phi$, where $x'$ denotes the variables of $\mathbb{R}^N$ which are orthogonal to $e$. Moreover, $\phi$ is decreasing in its first variable.\par
{\rm{(iii)}} If $\Omega=\R^N$ and $A$, $q$, $f(\cdot,s)$ $($for each $s\in\R)$, $p^{\pm}$ are constant, then $u$ is a planar $($i.e. one-dimensional$)$ traveling front, in the sense that
$$u(t,x)=\phi(x\cdot e-ct),$$
where $\phi\ :\ \R\to(p^-,p^+)$ is decreasing and $\phi(\mp\infty)=p^{\pm}$.
\end{theo}

Notice that properties (\ref{unifgamma}) and (\ref{tau}) are automatically satisfied here --and property~(\ref{tau}) is actually satisfied for all $\tau>0$-- due to the periodicity of $\Omega$, the definition of $\Omega^{\pm}_t$ and assumption (\ref{Gamma}).\par
The constant $\gamma(e)$ in (\ref{gamma}) is by definition larger than or equal to $1$. It measures the asymptotic ratio of the geodesic and Euclidean distances along the direction $e$. If the domain~$\Omega$ is invariant in the direction $e$, that is $\Omega=\Omega+se$ for all $s\in\R$, then $\gamma(e)=1$. For a pulsating traveling front satisfying (\ref{ptf}), the ``Euclidean speed" $c/\gamma(e)$ in the direction of propagation $e$ is then less than or equal to the global mean speed $c$ (the latter being indeed defined through the geodesic distance in $\Omega$).\par
Part (ii) of Theorem~\ref{cor1} still holds if $e$ is any direction of $\R^N$ and if $\Omega$, $A$, $q$, $f$, $\mu$ and $p^{\pm}$ are invariant in the direction~$e$ and periodic in the variables $x'$. This result can actually be extended to the case when the medium may not be periodic and $u$ may not be an invasion front:

\begin{theo}\label{statio}
Assume that $\Omega$ is invariant in a direction $e\in\mathbb{S}^{N-1}$, that $A$, $q$, $\mu$ and $p^{\pm}$ depend only on the variables $x'$ which are orthogonal to $e$, that $f=f(x',u)$ and that $(\ref{fdecreasing})$ and $(\ref{kappa})$ hold.\par
If $u$ is almost planar in the direction $e$, i.e. the sets $\Omega^{\pm}_t$ can be chosen as
$$\Omega^{\pm}_t=\big\{x\in\Omega;\ \pm(x\cdot e-\xi_t)<0\big\},$$
and if $u$ has global mean speed $c\ge 0$ with the stronger property that
$$\sup\big\{\,\big|\,|\xi_t-\xi_s|-c|t-s|\,\big|;\ (t,s)\in\R^2\big\}<+\infty,$$
then there exists $\epsilon\in\{-1,1\}$ such that
$$u(t,x)=\phi(x\cdot e-\epsilon ct,x')$$
for some function $\phi$. Moreover, $\phi$ is decreasing in its first variable.\par
If one further assumes that $c=0$, then the conclusion holds even if $f$ and $p^{\pm}$ also depend on $x\cdot e$, provided that they are nonincreasing in $x\cdot e$. In particular, if $u$ is quasi-stationary in the sense of Definition~$\ref{defspeed}$, then $u$ is stationary.
\end{theo}

In Theorems~\ref{cor1} and~\ref{statio}, we gave some conditions under which the fronts reduce to usual pulsating or traveling fronts. The fronts were assumed to have a global mean speed. Now, the following result generalizes part~(iii) of Theorem~\ref{cor1} to the case of almost planar fronts which may not have any global mean speed and which may not be invasion fronts. It gives some conditions under which almost planar fronts actually reduce to one-dimensional fronts.

\begin{theo}\label{th3} Assume that $\Omega=\R^N$, that $A$ and $q$ depend only on $t$, that the functions $p^{\pm}$ depend only on $t$ and $x\cdot e$ and are nonincreasing in $x\cdot e$ for some direction $e\in{\mathbb{S}}^{N-1}$, that $f=f(t,x\cdot e,u)$ is nonincreasing in $x\cdot e$, and that $(\ref{fdecreasing})$ and $(\ref{kappa})$ hold. If $u$ is almost planar in the direction $e$ with
$$\Omega^{\pm}_t=\big\{x\in\R^N;\ \pm(x\cdot e-\xi_t)<0\big\}$$
such that
\be\label{xisigma}
\exists\,\sigma>0,\ \sup\big\{|\xi_t-\xi_s|;\ (t,s)\in\R^2,\ |t-s|\le\sigma\big\}<+\infty,
\ee
then $u$ is planar, i.e. $u$ only depends on $t$ and $x\cdot e$~:
$$u(t,x)=\phi(t,x\cdot e)$$
for some function $\phi\ :\ \R^2\to\R$. Furthermore, 
\be\label{ineqsbis}
\forall\ (t,x)\in\R\times\R^N,\quad p^-(t,x\cdot e)<u(t,x)<p^+(t,x\cdot e)
\ee
and $u$ is decreasing with respect to $x\cdot e$.
\end{theo}

Notice that the assumption $\sup\,\{|\xi_{t+\sigma}-\xi_t|;\ t\in\R\}<+\infty$ for every $\sigma\in\R$ is clearly stronger than property (\ref{tau}). But one does not need $\xi_t$ to be monotone or $|\xi_t-\xi_s|\to+\infty$ as $|t-s|\to+\infty$, namely $u$ may not be an invasion front.\par
As for Theorem \ref{th1}, if the inequalities (\ref{ineqsbis}) are assumed to hold a priori and if $f$ is assumed to be nonincreasing in $s$ for $s$ in $[p^-(t,x\cdot e),p^-(t,x\cdot e)+\delta]$ and $[p^+(t,x\cdot e)-\delta,p^+(t,x\cdot e)]$ only, instead of $(-\infty,p^-(t,x\cdot e)+\delta]$ and $[p^+(t,x\cdot e)-\delta,+\infty)$, then the strict monotonicity of $u$ in $x\cdot e$ still holds.\par
As a particular case of the result stated in Theorem~\ref{statio} (with $c=0$), the following property holds, which states that, under some assumptions, any quasi-stationary front is actually stationary.

\begin{cor} Under the conditions of Theorem $\ref{th3}$, if one further assumes that the function $t\mapsto\xi_t$ is bounded and that $A$, $q$, $f$ and $p^{\pm}$ do not depend on $t$, then $u$ depends on~$x\cdot e$ only, that is $u$ is a stationary one-dimensional front.
\end{cor}


\subsection{Further extensions}\label{sec14}

In the previous sections, the waves were defined as spatial transitions connecting {\it two} limiting states $p^-$ and $p^+$. Multiple transition waves can be defined similarly.

\begin{defi}\label{multiple} {\rm{(Waves with multiple transitions)}} Let $k\ge 1$ be an  integer and let $p^1,\ldots,p^k$ be $k$ time-global solutions of $(\ref{eq})$. A {\rm{generalized transition wave connecting $p^1,\ldots,p^k$}} is a time-global classical solution $u$ of $(\ref{eq})$ such that $u\not\equiv p^j$ for all $1\le j\le k$, and there exist $k$ families $(\Omega^j_t)_{t\in{\mathbb R}}$ $(1\le j\le k)$ of open nonempty unbounded subsets of $\Omega$, a family $(\Gamma_t)_{t\in\R}$ of nonempty subsets of $\Omega$ and an integer $n\ge 1$ such that
$$\left\{\baa{l}
\forall\,t\in\R,\ \forall\,j\neq j'\in\{1,\ldots,k\},\ \Omega^j_t\cap\Omega^{j'}_t=\emptyset,\vspace{5pt}\\
\forall\,t\in\R,\ \displaystyle{\mathop{\bigcup}_{1\le j\le k}}(\partial\Omega^j_t\cap\Omega)=\Gamma_t,\ \ \Gamma_t\cup\displaystyle{\mathop{\bigcup}_{1\le j\le k}}\Omega^j_t=\Omega,\vspace{5pt}\\
\forall\,t\in\R,\ \forall\ j\in\{1,\ldots,k\},\ \sup\big\{d_{\Omega}(x,\Gamma_t);\ x\in\Omega^j_t\big\}=+\infty,\vspace{5pt}\\
\forall\,A\ge 0,\ \exists\,r>0,\ \forall\,t\in\R,\ \forall\,x\in\Gamma_t,\,\exists\,1\le j\neq j'\le k,\ \exists\,y^j\in\Omega^j_t,\ \exists\,y^{j'}\in\Omega^{j'}_t,\vspace{5pt}\\
\qquad\qquad\qquad\qquad\qquad d_{\Omega}(x,y^j)=d_{\Omega}(x,y^{j'})=r\hbox{ and }\min\big(d_{\Omega}(y^j,\Gamma_t),d_{\Omega}(y^{j'},\Gamma_t)\big)\ge A,\vspace{5pt}\\
\hbox{if }N=1\hbox{ then }\Gamma_t\hbox{ is made of at most }n\hbox{ points},\vspace{5pt}\\
\hbox{if }N\ge 2\hbox{ then }(\ref{omegapmbis})\hbox{ is satisfied},\eaa\right.$$
and
$$u(t,x)-p^j(t,x)\to 0\ \hbox{ uniformly in }t\in\R\hbox{ as }d_{\Omega}(x,\Gamma_t)\to+\infty\hbox{ and }x\in\overline{\Omega^j_t}$$
for all $1\le j\le k$.
\end{defi}

Triple or more general multiple transition waves are indeed known to exist in some reaction-diffusion problems (see e.g. \cite{bgs,fm}). The above definition also covers the case of multiple wave trains.\par
On the other hand, the spatially extended pulses, as defined in Definition~\ref{frontpulse} with $p^-(t)=p^+(t)$, correspond to the special case $k=1$, $p^1=p^{\pm}(t)$ and $\Omega^1_t=\Omega^-_t\cup\Omega^+_t$ in the above definition. We say that they are extended since, for each time $t$, the set $\Gamma_t$ is unbounded in general. The usual notion of localized pulses can be viewed as a particular case of Definition~\ref{multiple}.

\begin{defi} {\rm{(Localized pulses)}} In Definition~$\ref{multiple}$, if $k=1$ and if
$$\sup\big\{{\rm{diam}}_{\Omega}(\Gamma_t);\ t\in\R\big\}<+\infty,$$
then we say that $u$ is a {\rm{localized pulse}}.
\end{defi}

In all definitions of this paper, the time interval $\R$ can be replaced with any interval $I\subset\R$. However, when $I\neq\R$, the sets $\Omega^{\pm}_t$ or $\Omega^j_t$ are not required to be unbounded, but one only requires that
$$\lim_{t\to+\infty}\Big(\sup\big\{d_{\Omega}(x,\Gamma_t);\ x\in\Omega^{\pm}_t\}\Big)=+\infty\ \hbox{ or }\lim_{t\to+\infty}\Big(\sup\big\{d_{\Omega}(x,\Gamma_t);\ x\in\Omega^j_t\}\Big)=+\infty,$$
in the case of double or multiple transitions, if $I\supset[a,+\infty)$ (resp.
$$\lim_{t\to-\infty}\Big(\sup\big\{d_{\Omega}(x,\Gamma_t);\ x\in\Omega^{\pm}_t\}\Big)=+\infty\ \hbox{ or }\lim_{t\to-\infty}\Big(\sup\big\{d_{\Omega}(x,\Gamma_t);\ x\in\Omega^j_t\}\Big)=+\infty$$
if $I\supset(-\infty,a]$). The particular case $I=[0,T)$ with $0<T\le+\infty$ is used to describe the formation of waves and fronts for the solutions of Cauchy problems.\par
For instance, consider equation (\ref{homo}) for $t\ge 0$, with a function $f\in C^1([0,1])$ such that $f(0)=f(1)=0$, $f>0$ in $(0,1)$ and $f'(0)>0$. If $u_0$ is in $C_c(\R^N)$ and satisfies $0\le u_0\le 1$ with $u_0\not\equiv 0$ and if $u(t,x)$ denotes the solution of (\ref{homo}) with initial condition $u(0,\cdot)=u_0$, then $0\le u(t,x)\le 1$ for all $t\ge 0$ and $x\in\R^N$ and it follows easily from \cite{aw,j} that there exists a continuous increasing function $[0,+\infty)\ni t\mapsto r(t)>0$ such that $r(t)/t\to c^*>0$ as $t\to+\infty$ and
$$\left\{\baa{l}
\displaystyle{\mathop{\lim}_{A\to+\infty}}\Big(\inf\big\{u(t,x);\ t\ge 0,\ r(t)\ge A,\ 0\le|x|\le r(t)-A\big\}\Big)=1,\vspace{3pt}\\
\displaystyle{\mathop{\lim}_{A\to+\infty}}\Big(\sup\big\{u(t,x);\ t\ge 0,\ |x|\ge r(t)+A\big\}\Big)=0,\eaa\right.$$
where $c^*>0$ is the minimal speed of planar fronts $\varphi(x-ct)$ ranging in $[0,1]$ and connecting~$0$ and~$1$ for this equation (in other words, the mini\-mal speed $c^*$ of planar fronts is also the spreading speed of the solutions $u$ in all directions). If we define
$$\Omega^-_t=\big\{x\in\R^N;\ |x|<r(t)\big\},\ \Omega^+_t=\big\{x\in\R^N;\ |x|>r(t)\big\}\hbox{ and }\Gamma_t=\big\{x\in\R^N;\ |x|=r(t)\big\}$$
for all $t\ge 0$, then the function $u(t,x)$ can be viewed as a transition invasion wave connecting $p^-=1$ and $p^+=0$ in the time interval $[0,+\infty)$. We also refer to \cite{bhn3} for further definitions and properties of the spreading speeds of the solutions of the Cauchy problem $u_t=\Delta u+f(u)$ with compactly supported initial conditions, in arbitrary domains $\Omega$ and no-flux boundary conditions.\par
It is worth pointing out that, for the one-dimensional equation $u_t=u_{xx}+f(u)$ in $\R$ with $C^1([0,1],\R)$ functions $f$ such that $f(0)=f(1)=0$, $f(s)>0$ and $f'(s)\le f(s)/s$ on $(0,1)$, there are solutions $u:[0,+\infty)\times\R\to[0,1],\ (t,x)\mapsto u(t,x)$ such that
$$u(t,-\infty)=1,\ u(t,+\infty)=0\hbox{ for all }t\ge 0,\hbox{ and }\lim_{t\to+\infty}\|u_x(t,\cdot)\|_{L^{\infty}(\R)}=0,$$
see \cite{hr1}. At each time $t$, $u(t,\cdot)$ connects $1$ to $0$, but since the solutions become uniformly flatter and flatter as time runs, they are examples of solutions which are not generalized fronts connecting~$1$ and $0$.\hfill\break

\noindent{\bf{Time-dependent domains and other equations.}} We point out that all these general definitions can be adapted to the case when the domain
$$\Omega=\Omega_t$$
depends on time~$t$.\par
Lastly, the general definitions of transition waves which are given in this paper also hold for other types of evolution equations
$$F[t,x,u,Du,D^2u,\cdots]=0$$
which may not be of the parabolic type and which may be non local. Here $Du$ stands for the gradient of $u$ with respect to all variables $t$ and $x$.\hfill\break

\noindent{\bf{Outline of the paper.}} The following sections are devoted to proving all the results we have stated here. Section~\ref{intrinsic} is concerned with level set properties and the intrinsic character of the global mean speed. In Section~\ref{sec3}, we prove Theorem~\ref{thtime} on the existence of genera\-lized transition waves for the time-dependent equation~(\ref{eqkpp}). Section~\ref{sec4} deals with the proof of the general time-monotonicity result (Theorem~\ref{th1}). Section~\ref{sec5} is concerned with the proofs of Theorems~\ref{th2} and~\ref{cor1} on comparison of almost planar invasion fronts and reduction to pulsating fronts in periodic media. Lastly, in Section~\ref{sec6}, we prove the remaining Theorems~\ref{statio} and~\ref{th3} concerned with almost planar fronts in media which are invariant or monotone in the direction of propagation.


\SE{Intrinsic character of the interface localization and the global mean speed}\label{intrinsic}

Given a generalized transition wave $u$, we can view the set $\Gamma_t$ as the continuous interface of $u$ at time $t$. Of course this set is not uniquely defined, however, as we shall prove here, its localization in terms of (\ref{lambda}) and (\ref{infsup}) is intrinsic. Thus, this gives a meaning to the ``interface" in this continuous problem (even though it is not a free boundary). This section is divided into two parts, the first one dealing with the properties of the level sets and the second one with the intrinsic character of the global mean speed.


\subsection{Localization of the level sets: proof of Theorem~\ref{pro1}}

Heuristically, the fact that $u$ converges to two distinct constant states $p^{\pm}$ in $\Omega^{\pm}_t$ uniformly as $d_{\Omega}(x,\Gamma_t)\to+\infty$ will force any level set to stay at a finite distance from the interfaces~$\Gamma_t$, and the solution $u$ to stay away from $p^{\pm}$ in tubular neighborhoods of $\Gamma_t$.\par
More precisely, let us first prove part~1 of Theorem~\ref{pro1}. Formula~(\ref{lambda}) is almost imme\-diate. Indeed, assume that the conclusion does not hold for some $\lambda\in(p^-,p^+)$. Then there exists a sequence $(t_n,x_n)_{n\in\N}$ in $\R\times\overline{\Omega}$ such that
$$u(t_n,x_n)=\lambda\hbox{ for all }n\in\N\ \hbox{ and }\ d_{\Omega}(x_n,\Gamma_{t_n})\to+\infty\hbox{ as }n\to+\infty.$$
Up to extraction of some subsequence, two cases may occur: either $x_n\in\overline{\Omega^-_{t_n}}$ and then $u(t_n,x_n)\to p^-$ as $n\to+\infty$, or $x_n\in\overline{\Omega^+_{t_n}}$ and then $u(t_n,x_n)\to p^+$ as $n\to+\infty$. In both cases, one gets a contradiction with the fact that $u(t_n,x_n)=\lambda\in(p^-,p^+)$.\par
Assume now that property~(\ref{infsup}) does not hold for some $C\ge 0$. One may then assume that there exists a sequence $(t_n,x_n)_{n\in\N}$ of points in $\R\times\overline{\Omega}$ such that
\be\label{hyp}
d_{\Omega}(x_n,\Gamma_{t_n})\le C\hbox{ for all }n\in\N\ \hbox{ and }\ u(t_n,x_n)\to p^-\hbox{ as }n\to+\infty
\ee
(the case where $u(t_n,x_n)\to p^+$ could be treated similarly). Since $d_{\Omega}(x_n,\Gamma_{t_n})\le C$ for all $n$, it follows from (\ref{tau}) that there exists a sequence $(\tilde{x}_n)_{n\in\N}$ such that
$$\tilde{x}_n\in\Gamma_{t_n-\tau}\hbox{ for all }n\in\N\ \hbox{ and }\ \sup\big\{d_{\Omega}(x_n,\tilde{x}_n);\ n\in\N\big\}<+\infty.$$
On the other hand, from Definition \ref{def1}, there exists $d>0$ such that
$$\forall\ t\in\R,\ \forall\ y\in\overline{\Omega^+_t},\ d_{\Omega}(y,\Gamma_t)\ge d\ \Longrightarrow\ u(t,y)\ge\frac{p^-+p^+}{2}.$$
From (\ref{unifgamma}), there exists $r>0$ such that, for each $n\in\N$, there exists a point $y_n\in\overline{\Omega^+_{t_n-\tau}}$ satisfying
$$d_{\Omega}(\tilde{x}_n,y_n)=r\hbox{ and }d_{\Omega}(y_n,\Gamma_{t_n-\tau})\ge d.$$
Therefore,
\be\label{p-+}
\forall\ n\in\N,\ u(t_n-\tau,y_n)\ge\frac{p^-+p^+}{2}.
\ee
But the sequence $(d_{\Omega}(x_n,y_n))_{n\in\N}$ is bounded and the function $v=u-p^-$ is nonnegative and is a classical global solution of an equation of the type
$$v_t=\nabla_x\cdot(A(t,x)\nabla_x v)+q(t,x)\cdot\nabla_x v+b(t,x)v\hbox{ in }\R\times\overline{\Omega}$$
for some bounded function $b$, with $\mu(t,x)\cdot\nabla_xv(t,x)=0$ on $\partial\Omega$. Furthermore, the function~$v$ has bounded derivatives, from standard parabolic estimates. Since $v(t_n,x_n)\to 0$ as $n\to+\infty$ from (\ref{hyp}), one concludes from the linear estimates that $v(t_n-\tau,y_n)\to 0$ as $n\to+\infty$.\footnote{We use here the fact that, since the domain $\Omega$ is assumed to be globally smooth, as well as all coefficients $A$, $q$ and $\mu$ of (\ref{eq}) and (\ref{mu}), in the sense given in Section~\ref{intro}, then, for every positive real numbers $\delta,\,\rho,\,\sigma,\,M,\,B$ and $\eta>0$, there exists a positive real number $\epsilon=\epsilon(\delta,\rho,\sigma,M,B,\eta)>0$ such that, for any $t_0\in\R$, for any $C^1$ path $P:[0,1]\to\overline{\Omega}$ whose length is less than $\delta$, for any nonnegative classical supersolution~$\overline{u}$~of
$$\overline{u}_t\ge\nabla_x\cdot(A(t,x)\nabla_x\overline{u})+q(t,x)\cdot\nabla_x\overline{u}+b(t,x)\overline{u}$$
in the set
$$E=E_{t_0,P,\rho,\sigma}=[t_0,t_0+\rho]\!\times\!\big\{x\in\overline{\Omega};\ d_{\Omega}(x,P([0,1]))\le\rho\big\}\,\cup\,[t_0,t_0+\sigma]\!\times\! \overline{B_{\Omega}(P(0),\rho)},$$
satisfying~(\ref{mu}) on $\partial E\cap(\R\times\partial\Omega)$, $\|\nabla_x\overline{u}\|_{L^{\infty}(E)}\le M$, $\|b\|_{L^{\infty}(E)}\le B$ and
$$\max\!\big\{\overline{u}(t_0,P(s));\,s\in[0,1]\big\}\ge\eta,$$
then $\overline{u}(t_0+\sigma,P(0))\ge\epsilon$.} But
$$v(t_n-\tau,y_n)\ge\frac{p^+-p^-}{2}>0$$
from (\ref{p-+}). One has then reached a contradiction. This gives the desired conclusion~(\ref{infsup}).\par
To prove part 2 of Theorem~\ref{pro1}, assume now that~(\ref{lambda}) and~(\ref{infsup}) hold and that there is $d_0>0$ such that the sets
$$\big\{(t,x)\in\R\times\overline{\Omega};\ x\in\overline{\Omega^+_t},\ d_{\Omega}(x,\Gamma_t)\ge d\big\}$$
and
$$\big\{(t,x)\in\R\times\overline{\Omega};\ x\in\overline{\Omega^-_t},\ d_{\Omega}(x,\Gamma_t)\ge d\big\}$$
are connected for all $d\ge d_0$. Denote
$$m^-=\liminf_{x\in\overline{\Omega^-_t},\ d_{\Omega}(x,\Gamma_t)\to+\infty}u(t,x)\ \hbox{ and }\ M^-=\limsup_{x\in\overline{\Omega^-_t},\ d_{\Omega}(x,\Gamma_t)\to+\infty}u(t,x).$$
One has $p^-\le m^-\le M^-\le p^+$.\par
Call $\lambda=(m^-+M^-)/2$. Assume now that $m^-<M^-$. Then $\lambda\in(p^-,p^+)$ and, from~(\ref{lambda}), there exists $C_0\ge 0$ such that
$$d_{\Omega}(x,\Gamma_t)<C_0\ \hbox{ for all }(t,x)\in\R\times\overline{\Omega}\hbox{ such that }u(t,x)=\lambda.$$
Furthermore, there exist some times $t_1,t_2\in\R$ and some points $x_1,x_2$ with $x_i\in\overline{\Omega^-_{t_i}}$ such that $u(t_1,x_1)<\lambda<u(t_2,x_2)$ and $d_{\Omega}(x_i,\Gamma_{t_i})\ge\max(C_0,d_0)$ for $i=1,2$. Since the set
$$\big\{(t,x)\in\R\times\overline{\Omega};\ x\in\overline{\Omega^-_t},\ d_{\Omega}(x,\Gamma_t)\ge\max(C_0,d_0)\big\}$$
is connected and the function $u$ is continuous in $\R\times\overline{\Omega}$, there would then exist $t\in\R$ and $x\in\overline{\Omega^-_t}$ such that $d_{\Omega}(x,\Gamma_t)\ge\max(C_0,d_0)$ and $u(t,x)=\lambda$. But this is in contradiction with the choice of $C_0$.\par
Therefore, $p^-\le m^-=M^-\le p^+$ and
$$u(t,x)\to m^-\hbox{ uniformly as }d_{\Omega}(x,\Gamma_t)\to+\infty\hbox{ and }x\in\overline{\Omega^-_t}.$$
Similarly,
$$u(t,x)\to m^+\in[p^-,p^+]\hbox{ uniformly as }d_{\Omega}(x,\Gamma_t)\to+\infty\hbox{ and }x\in\overline{\Omega^+_t}.$$
If $\max(m^-,m^+)<p^+$, then there is $\epsilon>0$ and $C\ge 0$ such that $u(t,x)\le p^+-\epsilon$ for all $(t,x)$ with $d_{\Omega}(x,\Gamma_t)\ge C$. But
$$\sup\ \{u(t,x);\ d_{\Omega}(x,\Gamma_t)\le C\}<p^+$$
because of (\ref{infsup}). Therefore, $\sup\ \{u(t,x);\ (t,x)\in\R\times\overline{\Omega}\}<p^+$, which contradicts the fact that the range of $u$ is the whole interval $(p^-,p^+)$. As a consequence,
$$\max(m^-,m^+)=p^+.$$
Similarly, one can prove that $\min(m^-,m^+)=p^-$.\par
Eventually, either $m^-=p^-$ and $m^+=p^+$, or $m^-=p^+$ and $m^+=p^-$, which means that~$u$ is a transition wave connecting $p^-$ and $p^+$ (or $p^+$ and $p^-$). That completes the proof of Theorem~\ref{pro1}.\hfill$\Box$


\subsection{Uniqueness of the global mean speed for a given transition wave}

This section is devoted to the proof of the intrinsic character of the global mean speed, when it exists, of a generalized transition wave in the general vectorial case $m\ge 1$, when $p^+$ and~$p^-$ are separated from each other.\hfill\break

\noindent{\bf{Proof of Theorem \ref{pro2}.}} We make here all the assumptions of Theorem \ref{pro2} and we call
$$\tilde{\Gamma}_t=\partial\tilde{\Omega}^-_t\cap\Omega=\partial\tilde{\Omega}^+_t\cap\Omega$$
for all $t\in\R$. We first claim that there exists $C\ge 0$ such that
$$d_{\Omega}(x,\tilde{\Gamma}_t)\le C\hbox{ for all }t\in\R\hbox{ and } x\in\Gamma_t.$$
Assume not. Then there is a sequence $(t_n,x_n)_{n\in\N}$ in $\R\times\overline{\Omega}$ such that
$$x_n\in\Gamma_{t_n}\hbox{ for all }n\in\N\ \hbox{ and }\ d_{\Omega}(x_n,\tilde{\Gamma}_{t_n})\to+\infty\hbox{ as }n\to+\infty.$$
Up to extraction of some subsequence, one can assume that $x_n\in\overline{\tilde{\Omega}^-_{t_n}}$ (the case where $x_n\in\overline{\tilde{\Omega}^+_{t_n}}$ could be handled similarly). Call
$$\epsilon=\inf\big\{|p^-(t,x)-p^+(t,x)|;\ (t,x)\in{\mathbb R}\times\overline{\Omega}\big\}>0$$
and let $A\ge 0$ be such that
$$|u(t,z)-p^+(t,z)|\le\frac{\epsilon}{2}\hbox{ for all }(t,z)\in\R\times\overline{\Omega}\hbox{ with }d_{\Omega}(z,\Gamma_t)\ge A\hbox{ and }z\in\overline{\Omega^+_t}.$$
From the condition (\ref{unifgamma}), there exist $r>0$ and a sequence $(y_n)_{n\in\N}$ such that
$$y_n\in\overline{\Omega^+_{t_n}},\ \ d_{\Omega}(x_n,y_n)=r\ \hbox{ and }\ d_{\Omega}(y_n,\Gamma_{t_n})\ge A$$
for all $n\in\N$. Therefore,
$$d_{\Omega}(y_n,\tilde{\Gamma}_{t_n})\to+\infty\hbox{ as }n\to+\infty$$
and $y_n\in\overline{\tilde{\Omega}^-_{t_n}}$ for $n$ large enough. As a consequence,
$$u(t_n,y_n)-p^-(t_n,y_n)\to 0\ \hbox{ as }n\to+\infty.$$
On the other hand, $d_{\Omega}(y_n,\Gamma_{t_n})\ge A$ and $y_n\in\overline{\Omega^+_{t_n}}$, whence
$$|u(t_n,y_n)-p^+(t_n,y_n)|\le\frac{\epsilon}{2}$$
for all $n\in\N$. It follows that
$$\limsup_{n\to+\infty}|p^-(t_n,y_n)-p^+(t_n,y_n)|\le\frac{\epsilon}{2}.$$
This contradicts the definition of $\epsilon$.\par
Therefore, there exists $C\ge 0$ such that
\be\label{gamma'}
\forall\ t\in\R,\ \forall\ x\in\Gamma_t,\ d_{\Omega}(x,\tilde{\Gamma}_t)\le C.
\ee
Let now $(t,s)\in\R^2$ be any couple of real numbers and let $\eta>0$ be any positive number. There exists $(x,y)\in\Gamma_t\times\Gamma_s$ such that $d_{\Omega}(x,y)\le d_{\Omega}(\Gamma_t,\Gamma_s)+\eta$. From (\ref{gamma'}), there exists $(\tilde{x},\tilde{y})\in\tilde{\Gamma}_t\times\tilde{\Gamma}_s$ such that
$$d_{\Omega}(x,\tilde{x})\le C+\eta\ \hbox{ and }\ d_{\Omega}(y,\tilde{y})\le C+\eta.$$
Thus, $d_{\Omega}(\tilde{x},\tilde{y})\le d_{\Omega}(\Gamma_t,\Gamma_s)+2C+3\eta$ and
$$d_{\Omega}(\tilde{\Gamma}_t,\tilde{\Gamma}_s)\le d_{\Omega}(\Gamma_t,\Gamma_s)+2C+3\eta.$$
Since $\eta>0$ was arbitrary, one gets that $d_{\Omega}(\tilde{\Gamma}_t,\tilde{\Gamma}_s)\le d_{\Omega}(\Gamma_t,\Gamma_s)+2C$ for all $(t,s)\in\R^2$. Hence,
$$\limsup_{|t-s|\to+\infty}\frac{d_{\Omega}(\tilde{\Gamma}_t,\tilde{\Gamma}_s)}{|t-s|}\le\limsup_{|t-s|\to+\infty}\frac{d_{\Omega}(\Gamma_t,\Gamma_s)}{|t-s|}=c.$$
With similar arguments, by permuting the roles of the sets $\Omega^{\pm}_t$ and $\tilde{\Omega}^{\pm}_t$, one can prove that
$$d_{\Omega}(\Gamma_t,\Gamma_s)\le d_{\Omega}(\tilde{\Gamma}_t,\tilde{\Gamma}_s)+2\tilde{C}$$
for all $(t,s)\in\R^2$ and for some constant $\tilde{C}\ge 0$. Thus,
$$c=\liminf_{|t-s|\to+\infty}\frac{d_{\Omega}(\Gamma_t,\Gamma_s)}{|t-s|}\le\liminf_{|t-s|\to+\infty}\frac{d_{\Omega}(\tilde{\Gamma}_t,\tilde{\Gamma}_s)}{|t-s|}.$$
As a conclusion, the ratio $d_{\Omega}(\tilde{\Gamma}_t,\tilde{\Gamma}_s)/|t-s|$ converges as $|t-s|\to+\infty$, and its limit is equal to $c$. The proof of Theorem~\ref{pro2} is thereby complete.\hfill$\Box$


\SE{Generalized transition waves for a time-dependent equation}\label{sec3}

In this section, we construct explicit examples of generalized invasion transition fronts connecting $0$ and $1$ for the one-dimensional equation~(\ref{eqkpp}) under the assumption~(\ref{ftime}). Namely, we do the\hfill\break

\noindent{\bf{Proof of Theorem~\ref{thtime}.}} The strategy consists in starting from a classical traveling front with speed $c_1$ for the nonlinearity $f_1$, that is for times $t\in(-\infty,t_1]$, and then in letting it evolve and in proving that the solution eventually moves with speed $c_2$ at large times. The key point is to control the exponential decay of the solution when it approaches the state $0$, between times $t_1$ and $t_2$.\par
For the nonlinearity $f_1$, there exists a family of traveling fronts $\varphi_{1,c}(x-ct)$ of the equation
$$u_t=u_{xx}+f_1(u),$$
where $\varphi_{1,c}:\R\to(0,1)$ satisfies $\varphi_{1,c}(-\infty)=1$ and $\varphi_{1,c}(+\infty)=0$, for each speed $c\in[c^*_1,+\infty)$. The minimal speed $c^*_1$ satisfies $c^*_1\ge 2\sqrt{f'_1(0)}$, see \cite{aw,hr}. Each $\varphi_{1,c}$ is decreasing and unique up to shifts (one can normalize $\varphi_{1,c}$ is such a way that $\varphi_{1,c}(0)=1/2$). Furthermore, if $c>c^*_1$, then
$$\varphi_{1,c}(s)\sim A_{1,c}e^{-\lambda_{1,c}s}\ \hbox{ as }s\to+\infty,$$
where $A_{1,c}$ is a positive constant and $\lambda_{1,c}>0$ has been defined in (\ref{lambda1c}). If $c=c^*_1$ and $c^*_1>2\sqrt{f'_1(0)}$, then the same property holds. If $c=c^*_1$ and $c^*_1=2\sqrt{f'_1(0)}$, then
$$\varphi_{1,c}(s)\sim(A_{1,c}s+B_{1,c})\, e^{-\lambda_{1,c}s}\ \hbox{ as }s\to+\infty,$$
where $A_{1,c}\ge 0$, and $B_{1,c}>0$ if $A_{1,c}=0$, see \cite{aw}.\par
Let any speed $c_1\in[c^*_1,+\infty)$ be given, let $\xi$ be any real number (which is just a shift parameter) and let $u$ be the solution of (\ref{eqkpp}) such that
$$u(t,x)=\varphi_{1,c_1}(x-c_1t+\xi)\ \hbox{ for all }t\le t_1\hbox{ and }x\in\R.$$
Define
\be\label{xt1}
x_t=c_1t\ \hbox{ for all }t\le t_1.
\ee
The function $u$ satisfies
\be\label{tf1}\left\{\baa{ll}
u(t,x)\to 1 & \hbox{ as }x-x_t\to-\infty,\vspace{5pt}\\
u(t,x)\to 0 & \hbox{ as }x-x_t\to+\infty,\eaa\right.\hbox{ uniformly w.r.t. }t\le t_1.
\ee\par
Let us now study the behavior of $u$ on the time interval $[t_1,t_2]$ and next on the interval~$[t_2,+\infty)$. From the strong parabolic maximum principle, there holds $0<u(t,x)<1$ for all~$(t,x)\in\R^2$. For each $t\ge t_1$, the function~$u(t,\cdot)$ remains decreasing in~$\R$ since~$f$ does not depend on~$x$. Furthermore, from standard parabolic estimates, the function~$u$ satisfies the limiting conditions
\be\label{limloc}
u(t,-\infty)=1\ \hbox{ and }\ u(t,+\infty)=0\ \hbox{ locally in }t\in\R,
\ee
since $f(t,0)=f(t,1)=0$. Therefore, setting
\be\label{xt2}
x_t=x_{t_1}=c_1t_1\ \hbox{ for all }t\in(t_1,t_2],
\ee
one gets that 
\be\label{tf2}\left\{\baa{ll}
u(t,x)\to 1 & \hbox{ as }x-x_t\to-\infty,\vspace{5pt}\\
u(t,x)\to 0 & \hbox{ as }x-x_t\to+\infty,\eaa\right.\hbox{ uniformly w.r.t. }t\in(t_1,t_2].
\ee\par
Let $\epsilon$ be any positive real number in $(0,\lambda_{1,c_1})$. From the definition of $u$ and the above results, it follows that there exists a constant $C_{\epsilon}>0$, which also depends on $\xi$, $A_{1,c_1}$ and~$B_{1,c_1}$, such that
$$u(t_1,x)\le\min\left(C_{\epsilon}\,e^{-(\lambda_{1,c_1}-\epsilon)x},1\right)\ \hbox{ for all }x\in\R.$$
Let $M$ be the nonnegative real number defined by
$$M=\sup_{(t,s)\in[t_1,t_2]\times(0,1]}\frac{f(t,s)}{s}.$$
This quantity is finite since $f$ is of class $C^1$ and $f(t,0)=0$ for all $t$. Denote
$$\alpha=\lambda_{1,c_1}-\epsilon+\frac{M}{\lambda_{1,c_1}-\epsilon}>0$$
and
$$\overline{u}(t,x)=\min\left(C_{\epsilon}\,e^{-(\lambda_{1,c_1}-\epsilon)\,(x-\alpha(t-t_1))},1\right)\ \hbox{ for all }(t,x)\in[t_1,t_2]\times\R.$$
The function $\overline{u}$ is positive and it satisfies $u(t_1,\cdot)\le\overline{u}(t_1,\cdot)$ in $\R$. Furthermore, for all $(t,x)\in[t_1,t_2]\times\R$, if $\overline{u}(t,x)<1$, then
$$\baa{rcl}
\overline{u}_t(t,x)-\overline{u}_{xx}(t,x)-f(t,\overline{u}(t,x)) & \!\!\!\ge\!\!\! & \overline{u}_t(t,x)-\overline{u}_{xx}(t,x)-M\,\overline{u}(t,x)\vspace{5pt}\\
& \!\!\!\!=\!\! & C_{\epsilon}\left[\alpha\,(\lambda_{1,c_1}-\epsilon)-(\lambda_{1,c_1}-\epsilon)^2-M\right]\,e^{-(\lambda_{1,c_1}-\epsilon)\,(x-\alpha(t-t_1))}\vspace{5pt}\\
& \!\!\!=\!\!\! & 0\eaa$$
from the definitions of $M$ and $\alpha$. Thus, $\overline{u}$ is a supersolution of (\ref{eqkpp}) on the time interval $[t_1,t_2]$ and it is above $u$ at time $t_1$. Therefore,
\be\label{ineqsup}
u(t,x)\le\overline{u}(t,x)\le C_{\epsilon}\,e^{-(\lambda_{1,c_1}-\epsilon)\,(x-\alpha(t-t_1))}\ \hbox{ for all }(t,x)\in[t_1,t_2]\times\R
\ee
from the maximum principle.\par
On the other hand, from the behavior of $\varphi_{1,c_1}$ at $+\infty$, there exists a constant $C'_{\epsilon}>0$ such that
$$u(t_1,x)\ge\min\left(C'_{\epsilon}\,e^{-(\lambda_{1,c_1}+\epsilon)x},\frac{1}{2}\right)\ \hbox{ for all }x\in\R.$$
Let $\underline{u}$ the solution of the heat equation $\underline{u}_t=\underline{u}_{xx}$ for all $t\ge t_1$ and $x\in\R$, with value
$$\underline{u}(t_1,x)=\min\left(C'_{\epsilon}\,e^{-(\lambda_{1,c_1}+\epsilon)x},\frac{1}{2}\right)\ \hbox{ for all }x\in\R$$
at time $t_1$. Since $f\ge 0$, it follows from the maximum principle that
\be\label{ineqinf0}
u(t,x)\ge\underline{u}(t,x)\ \hbox{ for all }(t,x)\in[t_1,+\infty)\times\R.
\ee
But, for all $x\in\R$,
$$\underline{u}(t_2,x)=\int_{-\infty}^{+\infty}p(t_2-t_1,x-y)\,\underline{u}(t_1,y)\,dy\ge\,C'_{\epsilon}\int_{x_{\epsilon}}^{+\infty}p(t_2-t_1,x-y)\,e^{-(\lambda_{1,c_1}+\epsilon)y}\,dy,$$
where $x_{\epsilon}$ is the unique real number such that $C'_{\epsilon}\,e^{-(\lambda_{1,c_1}+\epsilon)x_{\epsilon}}=1/2$ and $p(\tau,z)=(4\pi\tau)^{-1/2}e^{-z^2/(4\tau)}$ is the heat kernel. Thus, for all $x\ge x_{\epsilon}+\sqrt{4(t_2-t_1)}$, there holds
\be\label{ineqinf}\baa{rcl}
\underline{u}(t_2,x) & \ge & \displaystyle{\frac{C'_{\epsilon}}{\sqrt{4\pi(t_2-t_1)}}}\displaystyle{\int_{x-\sqrt{4(t_2-t_1)}}^{x+\sqrt{4(t_2-t_1)}}}e^{-\frac{(x-y)^2}{4(t_2-t_1)}-(\lambda_{1,c_1}+\epsilon)y}dy\vspace{5pt}\\
& \ge & \displaystyle{\frac{2\,C'_{\epsilon}\,e^{-1-(\lambda_{1,c_1}+\epsilon)\sqrt{4(t_2-t_1)}}}{\sqrt{\pi}}}\times e^{-(\lambda_{1,c_1}+\epsilon)x}.\eaa
\ee\par
It follows from (\ref{ineqsup}), (\ref{ineqinf0}) and (\ref{ineqinf}) that, for all $\epsilon\in(0,\lambda_{1,c_1})$, there exist two positive constants $C_{\epsilon}^{\pm}$ and a real number $X_{\epsilon}$ such that
$$C_{\epsilon}^+e^{-(\lambda_{1,c_1}+\epsilon)x}\le u(t_2,x)\le C_{\epsilon}^-e^{-(\lambda_{1,c_1}-\epsilon)x}\ \hbox{ for all }x\in[X_{\epsilon},+\infty).$$
Remember also that $0<u(t_2,x)<1$ for all $x\in\R$, and that $u(t_2,-\infty)=1$. Since $f(t,s)=f_2(s)$ for all $t\ge t_2$ and $s\in[0,1]$, the classical front stability results (see e.g. \cite{l,u}) imply that
\be\label{phic2}
\sup_{x\in\R}\big|u(t,x)-\varphi_{2,c_2}(x-c_2t+m(t))\big|\to 0\ \hbox{ as }t\to+\infty,
\ee
where $m'(t)\to 0$ as $t\to+\infty$, and $c_2>0$ is given by (\ref{c2}). Here, $\varphi_{2,c_2}$ denotes the profile of the front traveling with speed $c_2$ for the equation $u_t=u_{xx}+f_2(u)$, such that $\varphi_{2,c_2}(-\infty)=1$ and $\varphi_{2,c_2}(+\infty)=0$. Therefore, there exists $t_3>t_2$ such that the map $t\mapsto c_2t-m(t)$ is increasing in $[t_3,+\infty)$, and $c_2t_3-m(t_3)\ge c_1t_1$. Define
\be\label{xt3}
x_t=\left\{\baa{ll}
c_1t_1 & \hbox{ if }t\in(t_2,t_3),\vspace{5pt}\\
c_2t-m(t) & \hbox{ if }t\in[t_3,+\infty).\eaa\right.
\ee
It follows from (\ref{limloc}) and (\ref{phic2}) that
\be\label{tf3}\left\{\baa{ll}
u(t,x)\to 1 & \hbox{ as }x-x_t\to-\infty,\vspace{5pt}\\
u(t,x)\to 0 & \hbox{ as }x-x_t\to+\infty,\eaa\right.\hbox{ uniformly w.r.t. }t\in(t_2,+\infty).
\ee\par
Eventually, setting
$$\Omega^{\pm}_t=\big\{x\in\R;\ \pm(x-x_t)<0\big\}$$
and $\Gamma_t=\{x_t\}$ for each $t\in\R$, where the real numbers $x_t$'s are defined in (\ref{xt1}), (\ref{xt2}) and (\ref{xt3}), one concludes from (\ref{tf1}), (\ref{tf2}) and (\ref{tf3}) that the function $u$ is a generalized transition front connecting $p^-=0$ and $p^+=1$. Furthermore, since the map $t\mapsto x_t$ is nondecreasing and $x_t-x_s\to+\infty$ as $t-s\to+\infty$, this transition front $u$ is an invasion of $0$ by $1$. The proof of Theorem~\ref{thtime} is thereby complete.\hfill$\Box$


\SE{Monotonicity properties}\label{sec4}

This section is devoted to the proof of the time-monotonicity properties, that is Theorem~\ref{th1}. This result has its own interest and it is also one of the key points in the subsequent uniqueness and classification results. The proof uses several comparison lemmata and some versions of the sliding method with respect to the time variable. Let us first show the following

\begin{pro}\label{lem1} Under the assumptions of Theorem $\ref{th1}$, one has
$$\forall\ (t,x)\in\R\times\overline{\Omega},\quad p^-(t,x)<u(t,x)<p^+(t,x).$$
\end{pro}

\noindent{\bf{Proof.}} We only prove the inequality $p^-(t,x)<u(t,x)$, the proof of the second inequality is similar. Remember that $u$ and $p^-$ are globally bounded. Assume now that
$$m:=\inf\big\{u(t,x)-p^-(t,x);\ (t,x)\in\R\times\overline{\Omega}\big\}<0.$$
Let $(t_n,x_n)_{n\in\N}$ be a sequence in $\R\times\overline{\Omega}$ such that
$$u(t_n,x_n)-p^-(t_n,x_n)\to m<0\hbox{ as }n\to+\infty.$$
Since $p^+(t,x)-p^-(t,x)\ge\kappa>0$ for all $(t,x)\in\R\times\overline{\Omega}$, it follows from Definition \ref{def1} that the sequence $(d_{\Omega}(x_n,\Gamma_{t_n}))_{n\in\N}$ is bounded. From assumption (\ref{tau}), there exists a sequence of points $(\tilde{x}_n)_{n\in\N}$ such that the sequence $(d_{\Omega}(x_n,\tilde{x}_n))_{n\in\N}$ is bounded and $\tilde{x}_n\in\Gamma_{t_n-\tau}$ for every $n\in\N$. From Definition \ref{def1}, there exists $d\ge 0$ such that
$$\forall\ t\in\R,\ \forall\ z\in\overline{\Omega^+_t},\ \  \big(d_{\Omega}(z,\Gamma_t)\ge d\big)\Longrightarrow\big(u(t,z)\ge p^+(t,z)-\kappa\big).$$
From the condition (\ref{unifgamma}), there exist $r>0$ and a sequence $(y_n)_{n\in\N}$ of points in $\overline{\Omega}$ such that
$$y_n\in\overline{\Omega^+_{t_n-\tau}},\ d_{\Omega}(y_n,\tilde{x}_n)=r\hbox{ and }d_{\Omega}(y_n,\Gamma_{t_n-\tau})\ge d\hbox{ for all }n\in\N.$$
One then gets that
\be\label{utnyn}
u(t_n-\tau_n,y_n)\ge p^+(t_n-\tau,y_n)-\kappa
\ee
for all $n\in\N$.\par
Call
$$v(t,x)=p^-(t,x)+m$$
and
$$w(t,x)=u(t,x)-v(t,x)=u(t,x)-p^-(t,x)-m\ge 0$$
for every $(t,x)\in\R\times\overline{\Omega}$. Since $p^-$ solves (\ref{eq}), since $f(t,x,\cdot)$ is nonincreasing in $(-\infty,p^-(t,x)+\delta]$ for each $(t,x)\in\R\times\overline{\Omega}$,\footnote{Here, we actually just use the fact that $f(t,x,\cdot)$ is nonincreasing in $(-\infty,p^-(t,x)]$ for each $(t,x)\in\R\times\overline{\Omega}$.} and since $m<0$, the function $v$ solves
$$v_t\le\nabla_x\cdot(A(x)\nabla_xv)+q(x)\cdot\nabla_xv+f(t,x,v)\hbox{ in }\R\times\overline{\Omega}$$
(remember that $A$ and $f$ do not depend on $t$, but this property is actually not used here). In other words, $v$ is a subsolution for (\ref{eq}). But $u$ solves (\ref{eq}) and $f(t,x,s)$ is locally Lipschitz-continuous in $s$ uniformly in $(t,x)\in\R\times\overline{\Omega}$. There exists then a bounded function $b$ such that
$$w_t\ge\nabla_x\cdot(A(x)\nabla_xv)+q(x)\cdot\nabla_xv+b(t,x)w\hbox{ in }\R\times\overline{\Omega}.$$
Lastly, $w$ satisfies $\mu\cdot\nabla_xw=0$ on $\R\times\partial\Omega$. Since the sequences $(d_{\Omega}(x_n,\tilde{x}_n))_{n\in\N}$ and $(d_{\Omega}(y_n,\tilde{x}_n))_{n\in\N}$ are bounded, the sequence $(d_{\Omega}(x_n,y_n))_{n\in\N}$ is bounded as well. Thus, since $w\ge 0$ in $\R\times\overline{\Omega}$ and $w(t_n,x_n)\to 0$ as $n\to+\infty$, one gets, as in the proof of part~1 of Theorem~\ref{pro1}, that $w(t_n-\tau,y_n)\to 0$ as $n\to+\infty$. But $w(t_n-\tau,y_n)$ satisfies
$$\baa{rcl}
w(t_n-\tau,y_n) & = & u(t_n-\tau,y_n)-p^-(t_n-\tau,y_n)-m\vspace{5pt}\\
& \ge & p^+(t_n-\tau,y_n)-\kappa-p^-(t_n-\tau,y_n)-m\ge-m>0\eaa$$
for all $n\in\N$ because of (\ref{utnyn}). One has then reached a contradiction.\par
As a conclusion, $m\ge 0$, whence
$$u(t,x)\ge p^-(t,x)\hbox{ for all }(t,x)\in\R\times\overline{\Omega}.$$\par
If $u(t_0,x_0)=p^-(t_0,x_0)$ for some $(t_0,x_0)\in\R\times\overline{\Omega}$, then the strong parabolic maximum principle and Hopf lemma imply that $u(t,x)=p^-(t,x)$ for all $x\in\overline{\Omega}$ and $t\le t_0$, and then for all $t\in\R$ by uniqueness of the Cauchy problem for (\ref{eq}). But this is impossible since $p^+-p^-\ge\kappa>0$ in~$\R\times\overline{\Omega}$ and $u(t,x)-p^+(t,x)\to 0$ uniformly as $x\in\overline{\Omega^+_t}$ and $d_{\Omega}(x,\Gamma_t)\to+\infty$ (notice actually that for each $t\in\R$, there are some points $z_n\in\overline{\Omega^+_t}$ such that $d_{\Omega}(z_n,\Gamma_t)\to+\infty$ as $n\to+\infty$, from (\ref{omegapm})).\par
As already underlined, the proof of the inequality $u<p^+$ is similar.\hfill$\Box$\break\par

Let us now turn to the\hfill\break

\noindent{\bf{Proof of Theorem \ref{th1}.}} In the hypotheses (\ref{fdecreasing}) and (\ref{kappa}), one can assume without loss of generality that $0<2\delta\le\kappa$, even if it means decreasing $\delta$. In what follows, for any $s\in\R$, we define $u^s$ in $\R\times\overline{\Omega}$ by
$$\forall\,(t,x)\in\R\times\overline{\Omega},\quad u^s(t,x)=u(t+s,x).$$
The general strategy is to prove that $u^s\ge u$ in $\R\times\overline{\Omega}$ for all $s>0$ large enough, and then for all $s\ge 0$ by sliding $u$ with respect to the time variable.\par
First, from Definition \ref{def1}, there exists $A>0$ such that
\be\label{A}
\forall\ (t,x)\in\R\times\overline{\Omega},\left\{\baa{l}
\big(x\in\overline{\Omega^-_t}\hbox{ and }d_{\Omega}(x,\Gamma_t)\ge A\big)\Longrightarrow\big(u(t,x)\le p^-(t,x)+\delta\big),\vspace{5pt}\\
\big(x\in\overline{\Omega^+_t}\hbox{ and }d_{\Omega}(x,\Gamma_t)\ge A\big)\Longrightarrow\Big(u(t,x)\ge p^+(t,x)-\displaystyle{\frac{\delta}{2}}\Big).\eaa\right.
\ee
Since $p^+$ invades $p^-$, there exists $s_0>0$ such that
$$\forall\,t\in\R,\ \forall\,s\ge s_0,\quad\Omega^+_{t+s}\supset\Omega^+_t\hbox{ and }d_{\Omega}(\Gamma_{t+s},\Gamma_t)\ge 2A.$$
Fix any $t\in\R$, $s\ge s_0$ and $x\in\overline{\Omega}$. If $x\in\overline{\Omega^+_t}$, then $x\in\overline{\Omega^+_{t+s}}$ and $d_{\Omega}(x,\Gamma_{t+s})\ge 2A$ since any continuous path from $x$ to $\overline{\Gamma_{t+s}}$ in $\overline{\Omega}$ meets $\overline{\Gamma_t}$. On the other hand, if $x\in\overline{\Omega^-_t}$ and $d_{\Omega}(x,\Gamma_t)\le A$, then $d_{\Omega}(x,\Gamma_{t+s})\ge A$ and $x\in\overline{\Omega^+_{t+s}}$. In both cases, one then has that
$$u^s(t,x)=u(t+s,x)\ge p^+(t+s,x)-\frac{\delta}{2}\ge p^+(t,x)-\delta$$
since $p^+$ is nondecreasing in time. To sum up,
\be\label{s0}\baa{rl}
\forall\,s\ge s_0,\ \ \forall\,(t,x)\in\R\times\overline{\Omega}, & \big(x\in\overline{\Omega^+_t}\big)\hbox{ or }\big(x\in\overline{\Omega^-_t}\hbox{ and }d_{\Omega}(x,\Gamma_t)\le A\big)\vspace{5pt}\\
& \qquad\Longrightarrow\big(u^s(t,x)=u(t+s,x)\ge p^+(t,x)-\delta\big).\eaa
\ee

\begin{lem}\label{lem2}
Call
$$\omega^-_A=\{(t,x)\in\R\times\overline{\Omega};\ x\in\overline{\Omega^-_t}\hbox{ and }d_{\Omega}(x,\Gamma_t)\ge A\}.$$
For all $s\ge s_0$, one has
$$u^s\ge u\ \hbox{ in }\omega^-_A.$$
\end{lem}

\noindent{\bf{Proof.}} Fix $s\ge s_0$ and define
$$\epsilon^*=\inf\big\{\epsilon>0;\ u^s\ge u-\epsilon\hbox{ in }\omega^-_A\big\}.$$
Since $u$ is bounded, $\epsilon^*$ is a well-defined nonnegative real number and one has
\be\label{useps}
u^s\ge u-\epsilon^*\hbox{ in }\omega^-_A.
\ee
One only has to prove that $\epsilon^*=0$.\par
Assume by contradiction that $\epsilon^*>0$. There exist then a sequence $(\epsilon_n)_{n\in\N}$ of positive real numbers and a sequence of points $(t_n,x_n)_{n\in\N}$ in $\omega^-_A$ such that
\be\label{epsilonn}
\epsilon_n\to\epsilon^*\hbox{ as }n\to+\infty\ \hbox{ and }\ u^s(t_n,x_n)<u(t_n,x_n)-\epsilon_n\hbox{ for all }n\in\N.
\ee\par
We first note that, when $x\in\overline{\Omega^-_t}$ and $d_{\Omega}(x,\Gamma_t)=A$, then $u(t,x)\le p^-(t,x)+\delta$ from~(\ref{A}), while $u^s(t,x)\ge p^+(t,x)-\delta$ from (\ref{s0}). Hence
\be\label{frontiere}\baa{rcl}
u^s(t,x)-u(t,x) & \ge & p^+(t,x)-p^-(t,x)-2\delta\vspace{5pt}\\
& \ge & \kappa-2\delta\ge0\ \hbox{ when }x\in\overline{\Omega^-_t}\hbox{ and }d_{\Omega}(x,\Gamma_t)=A.\eaa
\ee
Since $\nabla_xu$ is globally bounded in $\R\times\overline{\Omega}$, it follows from (\ref{epsilonn}) and the positivity of $\epsilon^*$ that there exists $\rho>0$ such that
$$\liminf_{n\to+\infty}d_{\Omega}(x_n,\Gamma_{t_n})\ge A+2\rho.$$
Even if it means decreasing $\rho$, one can also assume without loss of generality that
$$0<\rho<\tau,$$
where $\tau$ is given in (\ref{tau}), and that
\be\label{defrho}
\rho\times\big(\|(u-p^+)_t\|_{L^{\infty}(\R\times\overline{\Omega})}+\|\nabla_x(u-p^+)\|_{L^{\infty}(\R\times\overline{\Omega})}\big)\le\frac{\epsilon^*}{2}
\ee
since $u$ and $p^+$ have bounded derivatives.\par
Next, we claim that the sequence $(d_{\Omega}(x_n,\Gamma_{t_n}))_{n\in\N}$ is bounded. Otherwise, up to extraction of some subsequence, one has
$$d_{\Omega}(x_n,\Gamma_{t_n})\to+\infty\hbox{ and then }u(t_n,x_n)-p^-(t_n,x_n)\to 0\hbox{ as }n\to+\infty.$$
But, from Proposition \ref{lem1} and the fact that $p^-$ is nondecreasing in time, one has
$$\baa{rcl}
u(t_n,x_n)-p^-(t_n,x_n) & > & \epsilon_n+u(t_n+s,x_n)-p^-(t_n,x_n)\vspace{5pt}\\
& \ge & \epsilon_n+p^-(t_n+s,x_n)-p^-(t_n,x_n)\vspace{5pt}\\
& \ge & \epsilon_n\ \ \to\ \ \epsilon^*>0\hbox{ as }n\to+\infty,\eaa$$
which gives a contradiction. Therefore, the sequence $(d_{\Omega}(x_n,\Gamma_{t_n}))_{n\in\N}$ is bounded.\par
Since $x_n\in\overline{\Omega^-_{t_n}}$ and $d_{\Omega}(x_n,\Gamma_{t_n})\ge A+\rho$ for $n$ large enough (say, for $n\ge n_0$), and since~$p^+$ invades $p^-$, it follows that
$$x_n\in\overline{\Omega^-_t}\hbox{ and }d_{\Omega}(x_n,\Gamma_t)\ge A+\rho\ \hbox{ for all }n\ge n_0\hbox{ and }t\le t_n$$
and even that 
\be\label{inclusion}
x\in\overline{\Omega^-_t}\hbox{ and }d_{\Omega}(x,\Gamma_t)\ge A\ \hbox{ for all }n\ge n_0,\ x\in\overline{B_{\Omega}(x_n,\rho)}\hbox{ and }t\le t_n.
\ee
As a consequence, since $\rho<\tau$, there exists a sequence of points $(y_n)_{n\in\N,\,n\ge n_0}$ in $\overline{\Omega}$ such that
\be\label{yn}
y_n\in\overline{\Omega^-_{t_n-\tau+\rho}}\hbox{ and }A+\rho=d_{\Omega}(y_n,\Gamma_{t_n-\tau+\rho})=d_{\Omega}(x_n,\Gamma_{t_n-\tau+\rho})-d_{\Omega}(x_n,y_n)
\ee
for all $n\ge n_0$. Thus, for each $n\in\N$ with $n\ge n_0$, there exists a $C^1$ path $P_n:[0,1]\to\overline{\Omega^-_{t_n-\tau+\rho}}$ such that $P_n(0)=x_n$, $P_n(1)=y_n$, the length of $P_n$ is equal to $d_{\Omega}(x_n,y_n)$ and
$$d_{\Omega}(P_n(\sigma),\Gamma_{t_n-\tau+\rho})\ge A+\rho\hbox{ for all }\sigma\in[0,1].$$
Once again, since $p^+$ invades $p^-$, it follows that
\be\label{AP}
\forall\,n\ge n_0,\,\forall\,\sigma\in[0,1],\,\forall\,x\in\overline{B_{\Omega}(P_n(\sigma),\rho)},\,
\forall\,t\le t_n-\tau+\rho,\ x\in\overline{\Omega^-_t}\hbox{ and }d_{\Omega}(x,\Gamma_t)\ge A.
\ee
Together with (\ref{inclusion}), one gets that, for each $n\ge n_0$, the set
$$E_n=[t_n-\tau,t_n]\!\times\!\overline{B_{\Omega}(x_n,\rho)}\,\cup\,[t_n-\tau,t_n-\tau+\rho]\!\times\!\big\{x\in\overline{\Omega};\ d_{\Omega}(x,P_n([0,1]))\le\rho\big\}$$
is included in $\omega_A^-$.\par
As a consequence, for all $n\ge n_0$,
$$v:=u^s-(u-\epsilon^*)\ge 0\ \hbox{ in }E_n$$
from (\ref{useps}), and
$$u(t,x)-\epsilon^*<u(t,x)\le p^-(t,x)+\delta\ \hbox{ for all }(t,x)\in E_n$$
from (\ref{A}). Thus,
$$\baa{rcl}
(u-\epsilon^*)_t & = & \nabla_x\cdot(A(x)\nabla_x(u-\epsilon^*))+q(x)\cdot\nabla_x (u-\epsilon^*)+f(t,x,u)\vspace{5pt}\\
& \le & \nabla_x\cdot(A(x)\nabla_x(u-\epsilon^*))+q(x)\cdot\nabla_x (u-\epsilon^*)+f(t,x,u-\epsilon^*)\eaa$$
in $E_n$ for all $n\ge n_0$, because $f(t,x,\cdot)$ is nonincreasing in $(-\infty,p^-(t,x)+\delta]$. In other words, the function~$u-\epsilon^*$ is a subsolution of (\ref{eq}) in $E_n$ for all $n\ge n_0$. As far as the function $u^s(t,x)=u(t+s,x)$ is concerned, it satisfies
$$\baa{rcl}
u^s_t & = & \nabla_x\cdot(A(x)\nabla_xu^s)+q(x)\cdot\nabla_xu^s+f(t+s,x,u^s)\vspace{5pt}\\
& \ge & \nabla_x\cdot(A(x)\nabla_xu^s)+q(x)\cdot\nabla_xu^s+f(t,x,u^s)\eaa$$
for all $(t,x)\in\R\times\overline{\Omega}$ because $f(\cdot,x,\xi)$ is nondecreasing for all $(x,\xi)\in\overline{\Omega}\times\R$. Notice that we here use the fact that $A$ and $q$ are independent from the variable $t$. Furthermore, $u^s$ still satisfies
$$\mu(x)\cdot\nabla_xu^s(t,x)=0\hbox{ on }\R\times\partial\Omega$$
because $\mu$ is independent of $t$. In other words, $u^s$ is a supersolution of (\ref{eq}). Consequently, since the functions $f(t,x,\cdot)$ are locally Lipschitz-continuous uniformly with respect to $(t,x)\in\R\times\overline{\Omega}$, the function $v$ satisfies inequations of the type
$$v_t\ge\nabla_x\cdot(A(x)\nabla_xv)+q(x)\cdot\nabla_xv+b(t,x)v\ \hbox{ in }E_n$$
for all $n\ge n_0$, where the sequence $(\|b\|_{L^{\infty}(E_n)})_{n\in\N,\,n\ge n_0}$ is bounded.\par
On the other hand, since the sequence $(d_{\Omega}(x_n,\Gamma_{t_n}))_{n\in\N}$ is bounded, it follows from assumption (\ref{tau}) that there exists then a sequence of points $(\tilde{x}_n)_{n\in\N}$ in $\overline{\Omega}$ such that
$$\tilde{x}_n\in\Gamma_{t_n-\tau}\hbox{ for all }n\in\N,\hbox{ and }\sup\big\{d_{\Omega}(x_n,\tilde{x}_n);\ n\in\N\big\}<+\infty.$$
Thus, for all $n\ge n_0$,
$$d_{\Omega}(x_n,y_n)=d_{\Omega}(x_n,\Gamma_{t_n-\tau+\rho})-(A+\rho)\le d_{\Omega}(x_n,\Gamma_{t_n-\tau})-(A+\rho)\le d_{\Omega}(x_n,\tilde{x}_n)-(A+\rho)$$
since $x_n\in\overline{\Omega^-_{t_n}}$ and the sets $\Omega^-_t$ are non-increasing with respect to $t$ in the sense of the inclusion (because~$p^+$ invades $p^-$). The sequence $(d_{\Omega}(x_n,y_n))_{n\in\N,\,n\ge n_0}$ is then bounded. Lastly, remember that the function $\nabla_xv$ is bounded in $\R\times\overline{\Omega}$. As a conclusion, since $v(t_n,x_n)\to 0$ as $n\to+\infty$ (because of (\ref{epsilonn}) and $v(t_n,x_n)\ge 0$), it follows from the linear parabolic estimates that
\be\label{tnyn}
v(t_n-\tau,y_n)\to 0\ \hbox{ as }n\to+\infty.
\ee\par
But, because of (\ref{yn}), there exists a sequence $(z_n)_{n\in\N,\,n\ge n_0}$ such that
$$z_n\in\overline{\Omega^-_{t_n-\tau+\rho}},\ d_{\Omega}(y_n,z_n)=\rho\hbox{ and }d_{\Omega}(z_n,\Gamma_{t_n-\tau+\rho})=A$$
for all $n\ge n_0$. Thus, for all $n\ge n_0$,
$$u^s(t_n-\tau,y_n)-p^+(t_n-\tau,y_n)\ge u^s(t_n-\tau+\rho,z_n)-p^+(t_n-\tau+\rho,z_n)-\frac{\epsilon^*}{2}\ge-\delta-\frac{\epsilon^*}{2}$$
from (\ref{s0}) and (\ref{defrho}). Moreover,
$$u(t_n-\tau,y_n)\le p^-(t_n-\tau,y_n)+\delta\ \hbox{ for all }n\ge n_0$$
from (\ref{A}) and (\ref{AP}). Eventually, for all $n\ge n_0$, there holds
$$\baa{rcl}
v(t_n-\tau,y_n) & = & u^s(t_n-\tau,y_n)-u(t_n-\tau,y_n)+\epsilon^*\vspace{3pt}\\
& = & u^s(t_n-\tau,y_n)-p^+(t_n-\tau,y_n)+p^+(t_n-\tau,y_n)-u(t_n-\tau,y_n)+\epsilon^*\vspace{3pt}\\
& \ge & -\delta-\displaystyle{\frac{\epsilon^*}{2}}+p^+(t_n-\tau,y_n)-p^-(t_n-\tau,y_n)-\delta+\epsilon^*\vspace{3pt}\\
& \ge & \kappa-2\delta+\displaystyle{\frac{\epsilon^*}{2}}\ \ge\ \displaystyle{\frac{\epsilon^*}{2}}>0\eaa$$
from (\ref{kappa}) and the inequality $2\delta\le\kappa$.\par
One has then reached a contradiction with (\ref{tnyn}). Hence $\epsilon^*=0$ and the proof of Lemma~\ref{lem2} is thereby complete.\hfill$\Box$\break\par

Similarly, using now that $f(t,x,\cdot)$ is nonincreasing in $[p^+(t,x)-\delta,+\infty)$ and that $u^s(t,x)\ge p^+(t,x)-\delta/2\ge p^+(t,x)-\delta$ provided that $(t,x)\not\in\omega^-_A$ and $s\ge s_0$, we shall prove the following:

\begin{lem}\label{lem3}
For all $s\ge s_0$, one has
$$u^s\ge u\ \hbox{ in }\omega^+_A:=\R\times\overline{\Omega}\ \backslash\ \omega^-_A.$$
\end{lem}

\noindent{\bf{Proof.}} The proof uses some of the tools of that of Lemma~\ref{lem2}, but it is not just identical, because the time-sections of $\omega^+_A$, namely the sets $\overline{\Omega^+_t}\cup\big\{x\in\overline{\Omega^-_t};\ d_{\Omega}(x,\Gamma_t)<A\big\}$, are now nondecreasing with respect to time $t$ in the sense of the inclusion.\par
Fix $s\ge s_0$ and define
$$\epsilon_*=\inf\big\{\epsilon>0;\ u^s+\epsilon\ge u\hbox{ in }\omega^+_A\big\}.$$
This nonnegative real number is well-defined since $u$ is globally bounded, and one has
$$w:=u^s+\epsilon_*-u\ge 0\ \hbox{ in }\omega^+_A.$$
Furthermore, Lemma~\ref{lem2} implies that
\be\label{wepsilon}
w\ge\epsilon_*\hbox{ in }\omega^-_A.
\ee
In particular, $w$ is nonnegative in $\R\times\overline{\Omega}$.\par
To get the conclusion of Lemma~\ref{lem3}, it is sufficient to prove that $\epsilon_*=0$. Assume by contradiction that $\epsilon_*>0$. There exists then a sequence $(\epsilon_n)_{n\in\N}$ of positive real numbers and a sequence of points $(t_n,x_n)_{n\in\N}$ in $\omega^+_A$ such that
$$\epsilon_n\to\epsilon_*\hbox{ as }n\to+\infty,\hbox{ and }u^s(t_n,x_n)+\epsilon_n<u(t_n,x_n)\hbox{ for all }n\in\N.$$\par
If the sequence $(d_{\Omega}(x_n,\Gamma_{t_n}))_{n\in\N}$ were not bounded, then, up to extraction of a subsequence, it would converge to $+\infty$, whence
$$x_n\in\overline{\Omega^+_{t_n}}\subset\overline{\Omega^+_{t_n+s}}\hbox{ and }d_{\Omega}(x_n,\Gamma_{t_n+s})\ge d_{\Omega}(x_n,\Gamma_{t_n})\ \hbox{ for large }n.$$
Therefore, $d_{\Omega}(x_n,\Gamma_{t_n+s})\to+\infty$ and $u^s(t_n,x_n)-p^+(t_n+s,x_n)\to0$ as $n\to+\infty$. But
$$\baa{rcl}
u^s(t_n,x_n)-p^+(t_n+s,x_n) & < & u(t_n,x_n)-\epsilon_n-p^+(t_n+s,x_n)\vspace{5pt}\\
& \le & p^+(t_n,x_n)-p^+(t_n+s,x_n)-\epsilon_n\vspace{5pt}\\
& \le & -\epsilon_n\ \ \to\ \ -\epsilon_*<0\ \hbox{ as }n\to+\infty\eaa$$
from Proposition~\ref{lem1} and since $p^+$ is nondecreasing in time. This gives a contradiction.\par
Thus, the sequence $(d_{\Omega}(x_n,\Gamma_{t_n}))_{n\in\N}$ is bounded. From (\ref{tau}), there exists then a sequence~$(\tilde{x}_n)_{n\in\N}$ in $\overline{\Omega}$ such that
$$\tilde{x}_n\in\Gamma_{t_n-\tau}\hbox{ for all }n\in\N,\hbox{ and }\sup\big\{d_{\Omega}(x_n,\tilde{x}_n);\ n\in\N\big\}<+\infty.$$
Because of (\ref{unifgamma}), there exist $r>0$ and a sequence $(y_n)_{n\in\N}$ in $\overline{\Omega}$ such that
$$y_n\in\overline{\Omega^-_{t_n-\tau}},\ d_{\Omega}(\tilde{x}_n,y_n)=r\hbox{ and }d_{\Omega}(y_n,\Gamma_{t_n-\tau})\ge A\hbox{ for all }n\in\N.$$
There exists then a sequence $(z_n)_{n\in\N}$ in $\overline{\Omega}$ such that
\be\label{znA}
z_n\in\overline{\Omega^-_{t_n-\tau}}\hbox{ and }A=d_{\Omega}(z_n,\Gamma_{t_n-\tau})=d_{\Omega}(y_n,\Gamma_{t_n-\tau})-d_{\Omega}(y_n,z_n)\hbox{ for all }n\in\N.
\ee
Since $d_{\Omega}(y_n,z_n)\le d_{\Omega}(y_n,\Gamma_{t_n-\tau})\le d_{\Omega}(y_n,\tilde{x}_n)=r$ and since the sequence $(d_{\Omega}(x_n,\tilde{x}_n))_{n\in\N}$ is bounded, one gets finally that the sequence $(d_{\Omega}(x_n,z_n))_{n\in\N}$ is bounded.\par
Choose now $\rho>0$ so that
\be\label{defrho2}
\rho\,\|(u^s-u)_t\|_{L^{\infty}(\R\times\overline{\Omega})}+2\,\rho\,\|\nabla_x(u^s-u)\|_{L^{\infty}(\R\times\overline{\Omega})}<\epsilon_*
\ee
and $K\in\N\backslash\{0\}$ so that
\be\label{defK}
K\,\rho\ge\max\Big(\tau,\sup\big\{d_{\Omega}(x_n,z_n);\ n\in\N\big\}\Big).
\ee
For each $n\in\N$, there exists then a sequence of points $(X_{n,0},X_{n,1},\ldots,X_{n,K})$ in $\overline{\Omega}$ such that
$$X_{n,0}=x_n,\ X_{n,K}=z_n\hbox{ and }d_{\Omega}(X_{n,i},X_{n,i+1})\le\rho\hbox{ for each }0\le i\le K-1.$$
For each $n\in\N$ and $0\le i\le K-1$, set
$$E_{n,i}=\Big[t_n-\frac{i+1}{K}\,\tau, t_n-\frac{i}{K}\,\tau\Big]\times\overline{B_{\Omega}(X_{n,i},2\,\rho)}.$$\par
Since $w(t_n,x_n)\to 0$ as $n\to+\infty$, it follows from (\ref{defrho2}) and (\ref{defK}) that $w<\epsilon_*$ in $E_{n,0}$ for large $n$, whence $E_{n,0}\subset\omega^+_A$ from (\ref{wepsilon}). Consequently,
$$u^s(t,x)+\epsilon_*>u^s(t,x)\ge p^+(t,x)-\delta\hbox{ in }E_{n,0}\hbox{ for large }n$$
from (\ref{s0}). Since $f(t,x,\cdot)$ is nonincreasing in $[p^+(t,x)-\delta,+\infty)$ for all $(t,x)\in\R\times\overline{\Omega}$ and since $u^s$ is a supersolution of (\ref{eq}), it follows then as in the proof of Lemma~\ref{lem2} that the nonnegative function $w$ satisfies inequations of the type
$$w_t\ge\nabla_x(A(x)\nabla_xw)+q(x)\cdot\nabla_xw+b(t,x)w\ \hbox{ in }E_{n,0}$$
for $n$ large enough, where the sequence $(\|b\|_{L^{\infty}(E_{n,0})})_{n\in\N}$ is bounded. Remember also that $\mu(x)\cdot\nabla_xw(t,x)=0$ for all $(t,x)\in\R\times\partial\Omega$, and that $\nabla_xw$ is bounded in $\R\times\overline{\Omega}$. Since $w(t_n,X_{n,0})=w(t_n,x_n)\to 0$ as $n\to+\infty$, one concludes from the linear parabolic estimates that
$$w\Big(t_n-\frac{\tau}{K},X_{n,1}\Big)\to 0\hbox{ as }n\to+\infty.$$\par
An immediate induction yields $w(t_n-i\tau/K,X_{n,i})\to 0$ as $n\to+\infty$ for each $i=1,\ldots,K$. In particular, for $i=K$,
$$w(t_n-\tau,z_n)\to 0\hbox{ as }n\to+\infty.$$
But $z_n\in\overline{\Omega^-_{t_n-\tau}}$ and $d_{\Omega}(z_n,\Gamma_{t_n-\tau})=A$ for all $n\in\N$. As a consequence, for all $n\in\N$, $(t_n-\tau,z_n)\in\omega^-_A$ and $w(t_n-\tau,z_n)\ge\epsilon_*$ from (\ref{wepsilon}).\par
One has then reached a contradiction, which means that $\epsilon_*=0$. That completes the proof of Lemma~\ref{lem3}.\hfill$\Box$\break

\noindent{\bf{End of the proof of Theorem \ref{th1}.}} It follows from Lemmata \ref{lem2} and \ref{lem3} that
$$u^s\ge u\hbox{ in }\R\times\overline{\Omega}\hbox{ for all }s\ge s_0.$$
Now call 
$$s^*=\inf\ \{s>0;\ u^{\sigma}\ge u\hbox{ in }\R\times\overline{\Omega}\hbox{ for all }\sigma\ge s\}.$$
One has $0\le s^*\le s_0$ and one shall prove that $s^*=0$. Assume that $s^*>0$. Since $u^{s^*}\ge u$ in $\R\times\overline{\Omega}$, two cases may occur: either $\inf\big\{u^{s^*}(t,x)-u(t,x);\ d_{\Omega}(x,\Gamma_t)\le A\big\}>0$ or $\inf\big\{u^{s^*}(t,x)-u(t,x);\ d_{\Omega}(x,\Gamma_t)\le A\big\}=0$.\par
{\it{Case 1:}} assume that
$$\inf\big\{u^{s^*}(t,x)-u(t,x);\ d_{\Omega}(x,\Gamma_t)\le A\big\}>0.$$
Since $u_t$ is globally bounded, there exists $\eta_0\in(0,s^*)$ such that
\be\label{eta}
\forall\,\eta\in[0,\eta_0],\ \forall\,(t,x)\in\R\times\overline{\Omega},\ \ \big(d(x,\Gamma_t)\le A\big)\Longrightarrow\big(u^{s^*-\eta}(t,x)\ge u(t,x)\big).
\ee
For each $\eta\in[0,\eta_0]$, one then has $u^{s^*-\eta}(t,x)\ge u(t,x)$ for all $(t,x)\in\R\times\overline{\Omega}$ such that $x\in\overline{\Omega^-_t}$ and $d_{\Omega}(x,\Gamma_t)=A$, while $u(t,x)\le p^-(t,x)+\delta$ if $x\in\overline{\Omega^-_t}$ and $d_{\Omega}(x,\Gamma_t)\ge A$ (i.e. $(t,x)\in\omega^-_A$) from (\ref{A}). Therefore, the same arguments as in Lemma \ref{lem2} imply that
\be\label{eta0}
\forall\ \eta\in[0,\eta_0],\ u^{s^*-\eta}\ge u\hbox{ in }\omega^-_A.
\ee\par
On the other hand,
$$\big(x\in\overline{\Omega^+_t}\hbox{ and }d_{\Omega}(x,\Gamma_t)\ge A\big)\Longrightarrow\Big(u^{s^*}(t,x)\ge u(t,x)\ge p^+(t,x)-\frac{\delta}{2}\Big)$$
from (\ref{A}). Hence, even if it means decreasing $\eta_0>0$, one can assume without loss of generality that
$$\forall\eta\in[0,\eta_0],\ \ \big(x\in\overline{\Omega^+_t}\hbox{ and }d_{\Omega}(x,\Gamma_t)\ge A\big)\Longrightarrow\big(u^{s^*-\eta}(t,x)\ge p^+(t,x)-\delta\big).$$
Notice that this is the place where we use the choice of $\delta/2\ (<\delta)$ in the second property of~(\ref{A}). Furthermore, remember from (\ref{eta}) and (\ref{eta0}) that, for all $\eta\in[0,\eta_0]$, $u^{s^*-\eta}(t,x)\ge u(t,x)$ for all $(t,x)\in\R\times\overline{\Omega}$ such that 
$x\in\overline{\Omega^-_t}$, or $x\in\overline{\Omega^+_t}$ and $d_{\Omega}(x,\Gamma_t)\le A$. As in Lemma~\ref{lem3}, one then gets that
$$\forall\eta\in[0,\eta_0],\ \ \big(x\in\overline{\Omega^+_t}\hbox{ and }d_{\Omega}(x,\Gamma_t)\ge A\big)\Longrightarrow\big(u^{s^*-\eta}(t,x)\ge u(t,x)\big).$$\par
One concludes that $u^{s^*-\eta}\ge u$ in $\R\times\overline{\Omega}$ for all $\eta\in[0,\eta_0]$. That contradicts the minimality of $s^*$ and case 1 is then ruled out.\par
{\it{Case 2:}} assume that
$$\inf\big\{u^{s^*}(t,x)-u(t,x);\ d_{\Omega}(x,\Gamma_t)\le A\big\}=0.$$
There exists then a sequence $(t_n,x_n)_{n\in\N}$ in $\R\times\overline{\Omega}$ such that
$$d_{\Omega}(x_n,\Gamma_{t_n})\le A\hbox{ and }u^{s^*}(t_n,x_n)-u(t_n,x_n)\to 0\hbox{ as }n\to+\infty.$$
Since $u^{s^*}$ is a supersolution of (\ref{eq}) in $\R\times\overline{\Omega}$ (as already noticed in the proof of Lemma \ref{lem2}) and since $u^{s^*}\ge u$ in $\R\times\overline{\Omega}$, it follows from the linear parabolic estimates that
$$u(t_n,x_n)-u(t_n-s^*,x_n)=u^{s^*}(t_n-s^*,x_n)-u(t_n-s^*,x_n)\to 0\hbox{ as }n\to+\infty.$$
By immediate induction, one has that
\be\label{k}
u(t_n,x_n)-u(t_n-ks^*,x_n)\to 0\hbox{ as }n\to+\infty
\ee
for each $k\in\N$.\par
Fix any $\epsilon>0$. Let $B_{\epsilon}>0$ be such that
$$\forall\,(t,x)\in\R\times\overline{\Omega},\ \ \big(x\in\overline{\Omega^-_t}\hbox{ and }d_{\Omega}(x,\Gamma_t)\ge B_{\epsilon}\big)\Longrightarrow\big(u(t,x)\le p^-(t,x)+\epsilon\big).$$
On the other hand, since $p^+$ invades $p^-$ and since the sequence $(d_{\Omega}(x_n,\Gamma_{t_n}))_{n\in\N}$ is bounded, there exists $m\in\N$ such that
$$x_n\in\overline{\Omega^-_{t_n-ms^*}}\hbox{ and }d_{\Omega}(x_n,\Gamma_{t_n-ms^*})\ge B_{\epsilon}\hbox{ for all }n\in\N.$$
Hence,
$$u(t_n-ms^*,x_n)\le p^-(t_n-ms^*,x_n)+\epsilon\le p^-(t_n,x_n)+\epsilon\hbox{ for all }n\in\N$$
since $p^-$ is nondecreasing in time. Together with (\ref{k}) applied to $k=m$, one concludes that
$$\limsup_{n\to+\infty}\left(u(t_n,x_n)-p^-(t_n,x_n)\right)\le\epsilon.$$
But $u\ge p^-$ from Proposition~\ref{lem1}, and $\epsilon>0$ was arbitrary. One obtains that
\be\label{tnxn}
u(t_n,x_n)-p^-(t_n,x_n)\to 0\hbox{ as }n\to+\infty.
\ee\par
Let now $B>0$ be such that
$$\forall\,(t,x)\in\R\times\overline{\Omega},\ \ \big(x\in\overline{\Omega^+_t}\hbox{ and }d_{\Omega}(x,\Gamma_t)\ge B\big)\Longrightarrow\Big(u(t,x)\ge p^+(t,x)-\frac{\kappa}{2}\Big),$$
where $\kappa>0$ has been defined in (\ref{kappa}). From assumption (\ref{tau}), and since the sequence $(d_{\Omega}(x_n,\Gamma_{t_n}))_{n\in\N}$ is bounded, there exists a sequence $(\tilde{x}_n)_{n\in\N}$ in $\overline{\Omega}$ such that
$$\tilde{x}_n\in\Gamma_{t_n-\tau}\hbox{ for all }n\in\N,\hbox{ and }\sup\big\{d_{\Omega}(x_n,\tilde{x}_n);\ n\in\N\big\}<+\infty.$$
Because of (\ref{unifgamma}), there exist $r>0$ and a sequence $(y_n)_{n\in\N}$ in $\overline{\Omega}$ such that
$$y_n\in\overline{\Omega^+_{t_n-\tau}},\ d_{\Omega}(y_n,\tilde{x_n})=r\hbox{ and }d_{\Omega}(y_n,\Gamma_{t_n-\tau})\ge B\hbox{ for all }n\in\N.$$
Thus,
$$u(t_n-\tau,y_n)\ge p^+(t_n-\tau,y_n)-\frac{\kappa}{2}\hbox{ for all }n\in\N.$$
Remember now that both $u\ge p^-$ are two bounded solutions of (\ref{eq}) and that $f(t,x,\xi)$ is locally Lipschitz-continuous in $\xi$, uniformly with respect to $(t,x)\in\R\times\overline{\Omega}$. Notice also that the sequence $(d_{\Omega}(x_n,y_n))_{n\in\N}$ is bounded. Since $u(t_n,x_n)-p^-(t_n,x_n)\to 0$ as $n\to+\infty$ because of (\ref{tnxn}), one concludes that
$$u(t_n-\tau,y_n)-p^-(t_n-\tau,y_n)\to 0\hbox{ as }n\to+\infty.$$
But
$$u(t_n-\tau,y_n)-p^-(t_n-\tau,y_n)\ge p^+(t_n-\tau,y_n)-\frac{\kappa}{2}-p^-(t_n-\tau,y_n)\ge\frac{\kappa}{2}>0$$
owing to the definition of $\kappa$. One has then reached a contradiction and case 2 is then ruled out too.\par
As a consequence, $s^*=0$ and
$$u^s\ge u\hbox{ in }\R\times\overline{\Omega}\hbox{ for all }s\ge 0.$$\par
Let us now prove that the inequality is strict if $s>0$. Choose any $s>0$ and assume that
$$u^s(t_0,x_0)=u(t_0,x_0)\hbox{ for some }(t_0,x_0)\in\R\times\overline{\Omega}.$$
Since $u^s\ (\ge u)$ is a supersolution of (\ref{eq}), one gets that
$$u^s(t,x)=u(t,x)\hbox{ for all }t\le t_0\hbox{ and }x\in\overline{\Omega}$$
from the strong parabolic maximum principle and Hopf lemma. Fix any $t\le t_0$ and $x\in\overline{\Omega}$. For all $k\in\N$, one then has
$$0\le u(t,x)-p^-(t,x)=u(t-ks)-p^-(t,x)\le u(t-ks,x)-p^-(t-ks,x)$$
because $p^-$ is nondecreasing in time. But the right-hand side converges to $0$ as $k\to+\infty$, because $s>0$ and because of Definition \ref{invasion} (here, $p^+$ invades $p^-$). It follows that $u(t,x)=p^-(t,x)$ for all $t\le t_0$ and $x\in\overline{\Omega}$, which is impossible because of Proposition~\ref{lem1}.\par
As a conclusion, $u^s(t,x)>u(t,x)$ for all $(t,x)\in\R\times\overline{\Omega}$ and $s>0$. That completes the proof of Theorem \ref{th1}.\hfill$\Box$


\SE{Uniqueness of the mean speed, comparison of almost planar fronts and reduction to pulsating fronts}\label{sec5}

In this section, we prove, under some appropriate assumptions, the uniqueness of the speed among all almost-planar invasion fronts, and that the transition fronts reduce in some standard situations to the usual planar or pulsating fronts. Let us first process with the\hfill\break

\noindent{\bf{Proof of Theorem \ref{th2}.}} Notice first that $c$ and $\tilde{c}$ are (strictly) positive. Indeed,
$$d_{\Omega}(\Gamma_t,\Gamma_s),\ \  d_{\Omega}(\tilde{\Gamma}_t,\tilde{\Gamma}_s)\to+\infty\ \hbox{ as }|t-s|\to+\infty,$$
and the quantities $d_{\Omega}(\Gamma_t,\Gamma_s)-c\,|t-s|$ and $d_{\Omega}(\tilde{\Gamma}_t,\tilde{\Gamma}_s)-\tilde{c}\,|t-s|$ are assumed to be bounded uniformly with respect to $(t,s)\in\R^2$.\par
One shall prove that $c=\tilde{c}$ and that $\tilde{u}$ is above $u$ up to shift in time. Assume that $\tilde{c}<c$ (the other case can be treated similarly by permuting the roles of $u$ and $\tilde{u}$). Define
$$v(t,x)=\tilde{u}\left(\frac{c}{\tilde{c}}\,t,x\right)$$
and notice that
$$v_t(t,x)=\frac{c}{\tilde{c}}\ \tilde{u}_t\left(\frac{c}{\tilde{c}}\,t,x\right)\ge\tilde{u}_t\left(\frac{c}{\tilde{c}}\,t,x\right)=\nabla_x\cdot(A(x)\nabla_xv(t,x))+q(x)\cdot\nabla_xv(t,x)+f(x,v(t,x))$$
because $c/\tilde{c}\ge 1$ and $\tilde{u}_t\ge 0$ from Theorem \ref{th1}. We also use the fact that both $A$, $q$ and $f$ are independent of $t$. Furthermore, $\mu(x)\cdot\nabla_xv(t,x)=0$ on $\R\times\partial\Omega$. Therefore, the function~$v$, as well as all its time-shifts, is a supersolution for (\ref{eq}). It also follows from Definition \ref{def1} that
\be\label{vppm}
v(t,x)-p^{\pm}(x)\to 0\hbox{ uniformly as }x\in\overline{\tilde{\Omega}^{\pm}_{ct/\tilde{c}}}\hbox{ and }d_{\Omega}(x,\Gamma_{ct/\tilde{c}})\to+\infty,
\ee
where
$$\tilde{\Omega}^{\pm}_{ct/\tilde{c}}=\{x\in\Omega,\ \pm(x\cdot e-\tilde{\xi}_{ct/\tilde{c}})<0\}\hbox{ and }\tilde{\Gamma}^{\pm}_{ct/\tilde{c}}=\{x\in\Omega,\ x\cdot e=\tilde{\xi}_{ct/\tilde{c}}\}.$$
Remember that the quantities
$$d_{\Omega}(\tilde{\Gamma}_{ct/\tilde{c}},\tilde{\Gamma}_{cs/\tilde{c}})-\tilde{c}\left|\frac{c}{\tilde{c}}\,t-\frac{c}{\tilde{c}}\,s\right|=d_{\Omega}(\tilde{\Gamma}_{ct/\tilde{c}},\tilde{\Gamma}_{cs/\tilde{c}})-c|t-s|$$
are bounded independently of $(t,s)\in\R^2$. As a consequence, the map
$$t\mapsto d_{\Omega}(\Gamma_t,\Gamma_0)-d_{\Omega}(\tilde{\Gamma}_{ct/\tilde{c}},\tilde{\Gamma}_0)$$
is bounded in $\R$. Furthermore, both $u$ and $\tilde{u}$ are almost planar invasion fronts ($p^+$ invades $p^-$) in the same direction $e$, whence the maps $t\mapsto\xi_t$ and $t\mapsto\tilde{\xi}_t$ are nondecreasing. Eventually, one gets that
\be\label{sup}
\sup\big\{d_{\Omega}(\tilde{\Gamma}_{ct/\tilde{c}},\Gamma_t);\ t\in\R\big\}<+\infty\hbox{ and }\sup\big\{|\tilde{\xi}_{ct/\tilde{c}}-\xi_t|;\ t\in\R\big\}<+\infty.
\ee\par
On the other hand, Definition \ref{def1} applied to $u$ implies that there exists $A>0$ such that
\be\label{Abis}
\forall\,(t,x)\in\R\times\overline{\Omega},\ \left\{\baa{l}
\big(x\in\overline{\Omega^-_t}\hbox{ and }d_{\Omega}(x,\Gamma_t)\ge A\big)\Longrightarrow\big(u(t,x)\le p^-(x)+\delta\big)\vspace{5pt}\\
\big(x\in\overline{\Omega^+_t}\hbox{ and }d_{\Omega}(x,\Gamma_t)\ge A\big)\Longrightarrow\Big(u(t,x)\ge p^+(x)-\displaystyle{\frac{\delta}{2}}\Big).\eaa\right.
\ee
Since $u$ and $\tilde{u}$ are almost planar in the same direction $e$ and since $\tilde{u}$ is an invasion of $p^-$ by~$p^+$, properties (\ref{vppm}) and (\ref{sup}) yield the existence of $s_0>0$ such that, for all $s\ge s_0$ and for all $(t,x)\in\R\times\overline{\Omega}$,
$$\big(x\in\overline{\Omega^+_t}\big)\hbox{ or }\big(x\in\overline{\Omega^-_t}\hbox{ and }d_{\Omega}(x,\Gamma_t)\le A\big)\Longrightarrow\big(v^s(t,x)=v(t+s,x)\ge p^+(x)-\delta\big).$$
Choose any $s\ge s_0$. Since $p^-\le u,\ v\le p^+$ (from Proposition~\ref{lem1}) and
$$0<2\delta\le\kappa:=\inf\big\{p^+(t,x)-p^-(t,x);\ (t,x)\in\R\times\overline{\Omega}\big\}$$
(even if it means decreasing $\delta$ without loss of generality), the arguments used in Lemma~\ref{lem2} imply that
$$u(t,x)\le v^s(t,x)\hbox{ in }\omega^-_A,\hbox{ i.e. for all }x\in\overline{\Omega^-_t}\hbox{ such that }d_{\Omega}(x,\Gamma_t)\ge A.$$
Therefore, the arguments used in the proof of Lemma \ref{lem3} similarly imply that
$$u(t,x)\le v^s(t,x)\hbox{ for all }(t,x)\in\R\times\overline{\Omega}\ \backslash\ \omega^-_A.$$
Thus,
$$u\le v^s\hbox{ in }\R\times\overline{\Omega}\hbox{ for all }s\ge s_0.$$\par
Call now
$$s^*=\inf\ \{s\in\R;\ u\le v^s\hbox{ in }\R\times\overline{\Omega}\}.$$
One has $s^*\le s_0$ and $s^*>-\infty$ because $p^-(x)<u(t,x)<p^+(x)$ for all $(t,x)\in\R\times\overline{\Omega}$ (from Theorem \ref{th1}) and
$$v^s(0,x_0)=\tilde{u}\Big(\frac{c}{\tilde{c}}\,s,x_0\Big)\to p^-(x_0)\hbox{ as }s\to-\infty$$
for all $x_0\in\overline{\Omega}$ (see Definition \ref{invasion}). There holds
$$u\le v^{s^*}\hbox{ in }\R\times\overline{\Omega}.$$
In particular,
\be\label{vs*}
\big(x\in\overline{\Omega^+_t}\hbox{ and }d_{\Omega}(x,\Gamma_t)\ge A\big)\Longrightarrow\Big(v^{s^*}(t,x)\ge u(t,x)\ge p^+(x)-\frac{\delta}{2}\Big).
\ee\par
Assume now that
\be\label{hypcase1}
\inf\big\{v^{s^*}(t,x)-u(t,x);\ d_{\Omega}(x,\Gamma_t)\le A\big\}>0.
\ee
The same property then holds when $s^*$ is replaced with $s^*-\eta$ for any $\eta\in[0,\eta_0]$ and $\eta_0>0$ small enough, since $v_t$ (like $\tilde{u}_t$) is globally bounded. From (\ref{vs*}), one can assume that $\eta_0>0$ is small enough so that
$$\big(x\in\overline{\Omega^+_t}\hbox{ and }d_{\Omega}(x,\Gamma_t)\ge A\big)\Longrightarrow\big(v^{s^*-\eta}(t,x)\ge p^+(x)-\delta\big)$$
for all $\eta\in[0,\eta_0]$. The first property of (\ref{Abis}) implies, as in Lemma \ref{lem2}, that
$$v^{s^*-\eta}(t,x)\ge u(t,x)\hbox{ for all }\eta\in[0,\eta_0]\hbox{ and }(t,x)\in\R\times\overline{\Omega}\hbox{ with }x\in\overline{\Omega^-_t}\hbox{ and }d_{\Omega}(x,\Gamma_t)\ge A.$$
The above inequality then holds for all $(t,x)\in\R\times\overline{\Omega}$ such that $x\in\overline{\Omega^-_t}$, or $x\in\overline{\Omega^+_t}$ and $d_{\Omega}(x,\Gamma_t)\le A$. As in Lemma \ref{lem3}, one then gets that
$$v^{s^*-\eta}(t,x)\ge u(t,x)\hbox{ for all }\eta\in[0,\eta_0]\hbox{ and }(t,x)\in\R\times\overline{\Omega}\hbox{ with }x\in\overline{\Omega^+_t}\hbox{ and }d_{\Omega}(x,\Gamma_t)\ge A.$$
Eventually,
$$v^{s^*-\eta}\ge u\hbox{ in }\R\times\overline{\Omega}$$
for all $\eta\in[0,\eta_0]$. That contradicts the minimality of $s^*$ and assumption (\ref{hypcase1}) is false.\par
Therefore,
$$\inf\big\{v^{s^*}(t,x)-u(t,x);\ d_{\Omega}(x,\Gamma_t)\le A\big\}=0.$$
Then, there exists a sequence $(t_n,x_n)\in\R\times\overline{\Omega}$ such that $d_{\Omega}(x_n,\Gamma_{t_n})\le A$ for all $n\in\N$ and
$$v^{s^*}(t_n,x_n)-u(t_n,x_n)\to 0\hbox{ as }n\to+\infty.$$
Because of (\ref{tau}), there exists a sequence $(\tilde{x}_n)_{n\in\N}$ in $\overline{\Omega}$ such that
$$\tilde{x}_n\in\Gamma_{t_n-\tau}\hbox{ for all }n\in\N,\hbox{ and }\sup\big\{d_{\Omega}(x_n,\tilde{x}_n);\ n\in\N\big\}<+\infty.$$
Since $v^{s^*}$ is a supersolution of (\ref{eq}) and $v^{s^*}\ge u$ in $\R\times\overline{\Omega}$, it follows from the linear parabolic estimates that
$$\max\big\{v^{s^*}(t,x)-u(t,x);\ t_n-\tau-1\le t\le t_n-\tau,\ d_{\Omega}(x,\tilde{x}_n)\le 1\big\}\to 0\hbox{ as }n\to+\infty$$
and, since the functions $v^{s^*}_t$, $v^{s^*}_{x_i}$, $v^{s^*}_{x_ix_j}$, $u_t$, $u_{x_i}$ and $u_{x_ix_j}$ are globally H\"older continuous in~$\R\times\overline{\Omega}$ for all $1\le i,j\le N$, one gets that
$$\baa{rcl}
\big|v^{s^*}_t(t_n-\tau,\tilde{x}_n)-u_t(t_n-\tau,\tilde{x}_n)\big|+\big|v^{s^*}_{x_i}(t_n-\tau,\tilde{x}_n)-u_{x_i}(t_n-\tau,\tilde{x}_n)\big| & & \vspace{5pt}\\
+\big|v^{s^*}_{x_ix_j}(t_n-\tau,\tilde{x}_n)-u_{x_ix_j}(t_n-\tau,\tilde{x}_n)\big| & \to & 0\hbox{ as }n\to+\infty\eaa$$
for all $1\le i,j\le N$. But
$$\left\{\baa{rcl}
\displaystyle{\frac{\tilde{c}}{c}}\ v^{s^*}_t & = & \nabla_x\cdot(A(x)\nabla_xv^{s^*})+q(x)\cdot\nabla_xv^{s^*}+f(x,v^{s^*}),\vspace{5pt}\\
u_t & = & \nabla_x\cdot(A(x)\nabla_xu)+q(x)\cdot\nabla_xu+f(x,u).\eaa\right.$$
Therefore, $(\tilde{c}/c-1)\,u_t(t_n-\tau,\tilde{x}_n)\to 0$ as $n\to+\infty$, whence
$$u_t(t_n-\tau,\tilde{x}_n)\to 0\hbox{ as }n\to+\infty,$$
because $0<\tilde{c}<c$.\par
On the other hand, there exists $A'>0$ such that
$$\big(x\in\overline{\Omega^+_t}\hbox{ and }d_{\Omega}(x,\Gamma_t)\ge A'\big)\Longrightarrow\Big(u(t,x)\ge p^+(x)-\frac{\kappa}{3}\Big),$$
where $\kappa$ was defined in (\ref{kappa}). From (\ref{tau}), there exists a sequence $(y_n)_{n\in\N}$ in $\overline{\Omega}$ such that
$$y_n\in\Gamma_{t_n-2\tau}\hbox{ for all }n\in\N,\hbox{ and }\sup\big\{d_{\Omega}(\tilde{x}_n,y_n);\ n\in\N\big\}<+\infty.$$
Because of (\ref{unifgamma}), there exist $r>0$ and a sequence $(z_n)_{n\in\N}$ in $\overline{\Omega}$ such that
$$z_n\in\overline{\Omega^+_{t_n-2\tau}},\ d_{\Omega}(z_n,y_n)=r\hbox{ and }d_{\Omega}(z_n,\Gamma_{t_n-2\tau})\ge A'$$
for all $n\in\N$. Thus,
\be\label{zn}
u(t_n-2\tau,z_n)\ge p^+(z_n)-\frac{\kappa}{3}\hbox{ for all }n\in\N.
\ee
Since the sequence $(d_{\Omega}(z_n,\tilde{x}_n))_{n\in\N}$ is bounded, since $u_t(t_n-\tau,\tilde{x}_n)\to 0$ as $n\to+\infty$ and since the globally $C^1(\R\times\overline{\Omega})$ nonnegative function $u_t$ satisfies
$$(u_t)_t=\nabla_x\cdot(A(x)\nabla_xu_t)+q(x)\cdot\nabla_xu_t+f_t(t,x,u)\,u_t\ \hbox{ in }\R\times\overline{\Omega}$$
with $\|f_t(\cdot,\cdot,u(\cdot,\cdot))\|_{L^{\infty}(\R\times\overline{\Omega})}<+\infty$ and $\mu(x)\cdot\nabla_xu_t=0$ on $\R\times\partial\Omega$, the linear parabolic estimates imply that
$$u_t(t_n-2\tau,z_n)\to 0\hbox{ as }n\to+\infty.$$\par
Let now $\epsilon$ be any positive real number. Since the function $u_t$ is globally $C^1(\R\times\overline{\Omega})$, there exist $\sigma>0$ and $n_0\in\N$ such that
$$0\ \le\ \max\big\{u_t(t,z_n);\ t\in[t_n-2\tau-\sigma,t_n-2\tau]\big\}\ \le\ \epsilon\ \hbox{ for all }n\ge n_0.$$
Remember that $y_n\in\Gamma_{t_n-2\tau}$ and $d_{\Omega}(z_n,\Gamma_{t_n-2\tau})\le d_{\Omega}(z_n,y_n)=r$. Since $u$ is an invasion front of $p^-$ by $p^+$, there exists $\sigma'>0$ ($\sigma'$ is independent of $n$ and $\epsilon$) such that
\be\label{sigma'}
u(t_n-2\tau-\sigma',z_n)\le p^-(z_n)+\frac{\kappa}{3}\hbox{ for all }n\in\N.
\ee
Since $u_t(t_n-2\tau,z_n)\to 0$ as $n\to+\infty$ and $u_t\ge 0$ in $\R\times\overline{\Omega}$, it follows that, if $\sigma'\ge\sigma$, then
$$0\ \le\ \max\big\{u_t(t,z_n);\ t\in[t_n-2\tau-\sigma',t_n-2\tau-\sigma]\big\}\ \to\ 0\ \hbox{ as }n\to+\infty,$$
and then is less than $\epsilon$ for $n\ge n_1$ (for some $n_1\in\N$). Therefore, in both cases $\sigma'\ge\sigma$ or $\sigma'\le\sigma$, one has
$$0\ \le\ \max\big\{u_t(t,z_n);\ t\in[t_n-2\tau-\sigma',t_n-2\tau]\big\}\ \le\ \epsilon\ \hbox{ for all }n\ge\max(n_0,n_1).$$
Hence
$$u(t_n-2\tau-\sigma',z_n)\le u(t_n-2\tau,z_n)\le u(t_n-2\tau-\sigma',z_n)+\sigma'\epsilon$$
for $n$ large enough, and then
$$u(t_n-2\tau,z_n)-u(t_n-2\tau-\sigma',z_n)\to 0\hbox{ as }n\to+\infty$$
because $\epsilon>0$ was arbitrary and $\sigma'$ was independent of $\epsilon$. But
$$u(t_n-2\tau,z_n)-u(t_n-2\tau-\sigma',z_n)\ge p^+(z_n)-\frac{\kappa}{3}-p^-(z_n)-\frac{\kappa}{3}\ge\frac{\kappa}{3}>0\hbox{ for all }n\in\N$$
because of (\ref{zn}), (\ref{sigma'}) and of the definition of $\kappa$ in (\ref{kappa}). One has then reached a contradiction.\par
As a consequence,
$$\tilde{c}\ge c.$$
The other inequality follows by reversing the roles of $u$ and $\tilde{u}$. Thus, $c=\tilde{c}$.\par
The above arguments also imply that, for $u$ and $\tilde{u}$ as in Theorem \ref{th2}, there exists (the smallest) $T\in\R$ such that $\tilde{u}(t+T,x)\ge u(t,x)$ for all $(t,x)\in\R\times\overline{\Omega}$. The strong parabolic maximum principle and Hopf lemma imply that either the inequality is strict everywhere, or the two functions $u$ and $\tilde{u}^T$ are identically equal. That completes the proof of Theorem \ref{th2}.\hfill$\Box$\break

Let us now turn to the proof of the reduction of almost planar invasion fronts to pulsating fronts in periodic media.\hfill\break

\noindent{\bf{Proof of Theorem~\ref{cor1}.}} To prove part (i), fix $k\in L_1\Z\times\cdots\times L_N\Z$. By periodicity, the function
$$\tilde{u}(t,x)=u(t,x+k)$$
is a solution of (\ref{eq}). Furthermore, $\tilde{u}$, like $u$, satisfies all assumptions of Theorem \ref{th2}. Thus, there exists (the smallest) $T\in\R$ such that
\be\label{T}
\tilde{u}(t+T,x)=u(t+T,x+k)\ge u(t,x)\hbox{ for all }(t,x)\in\R\times\overline{\Omega}
\ee
and there exists a sequence of points $(t_n,x_n)_{n\in\N}$ in $\R\times\overline{\Omega}$ such that
\be\label{Tbis}
(d_{\Omega}(x_n,\Gamma_{t_n}))_{n\in\N}\hbox{ is bounded and }u(t_n+T,x_n+k)-u(t_n,x_n)\to 0\hbox{ as }n\to+\infty.
\ee
It shall then follow that
\be\label{liminf}
\liminf_{n\to+\infty}\ |u(t_n,x_n)-p^{\pm}(x_n)|>0.
\ee
Indeed, assume for instance that, up to extraction of some subsequence, $u(t_n,x_n)-p^-(x_n)\to 0$ as $n\to+\infty$ (the case $u(t_n,x_n)-p^+(x_n)\to 0$ as $n\to+\infty$ could be handled similarly). Then
$$\max\big\{u(t_n-\tau,y)-p^-(y);\ d_{\Omega}(y,x_n)\le C\big\}\to 0\hbox{ as }n\to+\infty$$
for any $C\ge 0$, from the linear parabolic estimates applied to the nonnegative function $u-p^-$ (remember that $\tau>0$ is given in (\ref{tau})). But there is a sequence $(y_n)_{n\in\N}$ in $\overline{\Omega}$ such that
$$(d_{\Omega}(y_n,x_n))_{n\in\N}\hbox{ is bounded },y_n\in\overline{\Omega^+_{t_n-\tau}}\hbox{ and }u(t_n-\tau,y_n)\ge p^+(y_n)-\frac{\kappa}{2}\hbox{ for all }n\in\N$$
(one uses the facts that the sequence $(d_{\Omega}(x_n,\Gamma_{t_n}))_{n\in\N}$ is bounded and that (\ref{unifgamma}) is automa\-ti\-cally satisfied by periodicity of $\Omega$). One then gets a contradiction as $n\to+\infty$. Thus, (\ref{liminf}) holds.\par
Write
$$x_n=x'_n+x''_n$$
for all $n\in\N$, where $x'_n\in L_1\Z\times\cdots\times L_N\Z$ and $x''_n\in[0,L_1]\times\cdots\times[0,L_N]\ \cap\ \overline{\Omega}$. Set
$$u_n(t,x)=u(t+t_n,x+x'_n)$$
for all $n\in\N$ and $(t,x)\in\R\times\overline{\Omega}$. The functions $u_n$ satisfy the same equation (\ref{eq}) with the same boundary conditions (\ref{mu}) as $u$, since the domain $\Omega$ is periodic and the coefficients $A$, $q$, $f$ and $\mu$ are periodic and independent of $t$. Up to extraction of a subsequence one can assume that $x''_n\to x_{\infty}\in\overline{\Omega}$ as $n\to+\infty$ and that, from standard parabolic estimates, $u_n(t,x)\to u_{\infty}(t,x)$ locally uniformly in $\R\times\overline{\Omega}$, where $u_{\infty}$ solves (\ref{eq}) and (\ref{mu}). Furthermore,
$$u_{\infty}(t+T,x+k)\ge u_{\infty}(t,x)\hbox{ for all }(t,x)\in\R\times\overline{\Omega}$$
from (\ref{T}), and
$$u_{\infty}(T,x_{\infty}+k)=u_{\infty}(0,x_{\infty})$$
from (\ref{Tbis}). It follows then from the strong maximum principle, Hopf lemma and the uniqueness of the solution of the Cauchy problem for (\ref{eq}) and (\ref{mu}), that
\be\label{uinfty}
u_{\infty}(t+T,x+k)=u_{\infty}(t,x)\hbox{ for all }(t,x)\in\R\times\overline{\Omega}.
\ee
Furthermore,
\be\label{neq}
u_{\infty}(0,x_{\infty})\neq p^{\pm}(x_{\infty})
\ee
from (\ref{liminf}).\par
On the other hand, as already noticed in the proof of Theorem \ref{th2}, the global mean speed $c$ is positive. Since the quantities $d_{\Omega}(\Gamma_t,\Gamma_s)-c|t-s|$ are bounded independently of $(t,s)\in\R^2$, since
$$\Omega^{\pm}_t=\{x\in\Omega,\ \pm(x\cdot e-\xi_t)<0\}$$
and since $t\mapsto\xi_t$ is nondecreasing (because $p^+$ invades $p^-$), it follows from the definition of $\gamma=\gamma(e)$ in (\ref{gamma}) that there exists $M\ge 0$ such that
\be\label{xits}
\big|\xi_t-c\,\gamma^{-1}\,t\big|\le M\hbox{ for all }t\in\R.
\ee\par
But, from Definition~\ref{def1}, since the geodesic distance is not smaller than the Euclidean distance, one has that
$$u_n(t,x)-p^{\pm}(x)=u(t+t_n,x+x'_n)-p^{\pm}(x+x'_n)\to 0\ \hbox{ as }(x+x'_n)\cdot e-\xi_{t+t_n}\to\mp\infty,$$
uniformly with respect to $n$ and $(t,x)$. Write
$$(x+x'_n)\cdot e-\xi_{t+t_n}=x\cdot e-\xi_t+x_n\cdot e-\xi_{t_n}-x''_n\cdot e+\xi_t+\xi_{t_n}-\xi_{t+t_n}.$$
The sequence $(x_n\cdot e-\xi_{t_n})_{n\in\N}$ is bounded because $(d_{\Omega}(x_n,\Gamma_{t_n}))_{n\in\N}$ is bounded. The quantities $\xi_t+\xi_{t_n}-\xi_{t+t_n}$ are bounded independently of $t$ and $n$ because of (\ref{xits}). Lastly, the sequence~$(x''_n\cdot e)_{n\in\N}$ is also bounded. Finally, one gets that
$$u_{\infty}(t,x)-p^{\pm}(x)\to 0\ \hbox{ as }x\cdot e-\xi_t\to\mp\infty$$
uniformly with respect to $(t,x)$.\par
Assume now, by contradiction, that $T>\gamma(k\cdot e)/c$ (one shall actually prove that $T=\gamma(k\cdot e)/c$). Since
$$(x_{\infty}+mk)\cdot e-\xi_{mT}\to\mp\infty\ \hbox{ as }m\in\Z\hbox{ and }m\to\pm\infty$$
because of our assumption and because of (\ref{xits}), it follows that
$$u_{\infty}(mT,x_{\infty}+mk)-p^{\pm}(x_{\infty}+mk)\to 0\ \hbox{ as }m\in\Z\hbox{ and }m\to\pm\infty.$$
But $p^{\pm}(x_{\infty}+mk)=p^{\pm}(x_{\infty})$ for all $m\in\Z$ by periodicity of $p^{\pm}$, and $u_{\infty}(mT,x_{\infty}+mk)=u_{\infty}(0,x_{\infty})$ for all $m\in\Z$ because of (\ref{uinfty}). One finally gets a contradiction with (\ref{neq}).\par
Therefore, the inequality $T>\gamma(k\cdot e)/c$ was impossible, whence $T\le\gamma(k\cdot e)/c$ and
$$u\left(t+\frac{\gamma\ k\cdot e}{c},x+k\right)\ge u(t,x)\hbox{ for all }(t,x)\in\R\times\overline{\Omega}.$$\par
Similarly, by fixing the function $u(t,x+k)$ and sliding $u(t,x)$ with respect to $t$, one can prove that
$$u\left(t-\frac{\gamma\ k\cdot e}{c},x\right)\ge u(t,x+k)\hbox{ for all }(t,x)\in\R\times\overline{\Omega}.$$\par
As a consequence,
\be\label{ptfbis}
u\left(t+\frac{\gamma\ k\cdot e}{c},x+k\right)=u(t,x)\hbox{ for all }(t,x)\in\R\times\overline{\Omega},
\ee
namely $u$ is a pulsating traveling front in the sense of (\ref{ptf}). Its global mean speed is equal to $c\,\gamma^{-1}$ in the sense of (\ref{ptf}), but it is equal to $c$ in the more intrinsic sense of Definition~\ref{defspeed}.\par
Let now $u$ and $v$ be two fronts satisfying all assumptions of part (i) of Theorem~\ref{cor1}. One shall prove that $u$ and $v$ are equal up to shift in time. From Theorem~\ref{th2}, there exists (the smallest) $T\in\R$ such that
$$v(t+T,x)\ge u(t,x)\hbox{ for all }(t,x)\in\R\times\overline{\Omega}$$
and there exists a sequence $(t_n,x_n)_{n\in\N}$ in $\R\times\overline{\Omega}$ such that
$$(d_{\Omega}(x_n,\Gamma_{t_n}))_{n\in\N}\hbox{ is bounded, and }v(t_n+T,x_n)-u(t_n,x_n)\to 0\hbox{ as }n\to+\infty.$$
Since both $u$ and $v$ satisfy (\ref{ptfbis}) for all $k\in L_1\Z\times\cdots\times L_N\Z$, one can assume without loss of generality that the sequence $(x_n)_{n\in\N}$ is bounded. But since the sequence $(d_{\Omega}(x_n,\Gamma_{t_n}))_{n\in\N}$ is itself bounded and since $u$ is an invasion front, the sequence $(t_n)_{n\in\N}$ is then bounded as well. Up to extraction of some subsequence, one can then assume that $(t_n,x_n)\to(\overline{t},\overline{x})\in\R\times\overline{\Omega}$, whence
$$v(\overline{t}+T,\overline{x})=u(\overline{t},\overline{x}).$$
The strong parabolic maximum principle and Hopf lemma then yield
$$v(t+T,x)=u(t,x)\hbox{ for all }(t,x)\in\R\times\overline{\Omega},$$
which completes the proof of part (i) of Theorem~\ref{cor1}.\hfill\break\par
To prove part (ii), assume, without loss of generality, that $e=e_1=(1,0,\ldots,0)$. Fix any $\sigma\in\R\backslash\{0\}$. The data $\Omega$, $A$, $q$, $f$, $\mu$ and $p^{\pm}$ are then periodic with respect to the positive vector $(|\sigma|,L_2,\ldots,L_N)$. Part (i) applied to $k=(\sigma,0,\ldots,0)$ then implies that
$$u\left(t+\frac{\gamma\,\sigma}{c},x\right)=u(t,x_1-\sigma,x_2,\ldots,x_N)$$
for all $(t,x)\in\R\times\overline{\Omega}$, where $\gamma=\gamma(e)=1$ since $\Omega$ is invariant in the direction $e$. Since this property holds for any $\sigma\in\R\backslash\{0\}$ (and also for $\sigma=0$ obviously), it follows that
$$u(t,x)=\phi(x_1-ct,x')\hbox{ for all }(t,x)\in\R\times\overline{\Omega},$$
where $x'=(x_2,\ldots,x_N)$ and the function $\phi\ :\ \overline{\Omega}\to\R$ is defined by
$$\phi(\zeta,x')=u\left(-\frac{\zeta}{c},0,x'\right)\hbox{ for all }(\zeta,x')\in\overline{\Omega}.$$
The function $\phi$ is then decreasing in $\zeta$ since $u$ is increasing in $t$ and $c>0$.\hfill\break\par
Lastly, part (iii) is a consequence of part (ii) and of Theorem \ref{th3}. Namely, part (ii) implies that $u$ depends only on $x\cdot e-ct$ and on the variables $x'$ which are orthogonal to~$e$, and Theorem \ref{th3} (its proof will be done in Section~\ref{sec6}) implies that $u$ does not depend on~$x'$.\footnote{Notice that, in this part (iii), one can assume without loss of generality that $\xi_t=c\,t$ for all $t\in\R$, because of Definition~\ref{def1}.} Therefore, 
$$u(t,x)=\phi(x\cdot e-ct)\hbox{ for all }(t,x)\in\R\times\R^N,$$
where the function $\phi\ :\ \R\to\R$ is defined by $\phi(\zeta)=u(-\zeta/c,0,\ldots,0)$ for all $\zeta\in\R$, is decreasing in $\R$ and satisfies $\phi(\mp\infty)=p^{\pm}$. The proof of Theorem~\ref{cor1} is now complete.\hfill$\Box$


\SE{The case of media which are invariant or monotone in the direction of propagation}\label{sec6}

In this section, we assume that the domain is invariant in a direction~$e$ and we prove that, under appropriate conditions on the coefficients of (\ref{eq}), the almost planar fronts, which may not be invasions, do not depend on the transverse variables or have a constant profile in the direction~$e$. We start with the\hfill\break

\noindent{\bf{Proof of Theorem \ref{th3}.}} Up to rotation of the frame, one can assume without loss of generality that $e=e_1=(1,0,\ldots,0)$. We shall then  prove that $u$ is decreasing in $x_1$ and that it does not depend on the variable $x'=(x_2,\ldots,x_N)$.\par
First, notice that the same arguments as in Proposition~\ref{lem1} yield the inequalities (\ref{ineqsbis}). The proof is even simpler here due to the facts that $\Gamma_t=\{x_1=\xi_t\}$ and that assumption~(\ref{xisigma}) is made.\par
Actually, because of (\ref{xisigma}) and Definition~\ref{def1}, one can assume without loss of generality in the sequel that the map $t\mapsto\xi_t$ is uniformly continuous in $\R$.\par
Fix any vector $\theta\in\R^{N-1}$ and call
$$v(t,x)=u(t,x_1,x'+\theta).$$
Since the coefficients of (\ref{eq}) are assumed to be independent of $x'$, the function $v$ is a solution of the same equation (\ref{eq}) as $u$, with the same choice of sets $(\Omega^{\pm}_t)_{t\in\R}$ and $(\Gamma_t)_{t\in\R}$. Let $A\ge 0$ be such that
\be\label{Ater}\forall\,(t,x)\in\R\times\R^N,\ \left\{\baa{lcl}
\big(x_1-\xi_t\ge A\big) & \Longrightarrow & \big(u(t,x)\le p^-(t,x_1)+\delta\big)\vspace{5pt}\\
\big(x_1-\xi_t\le -A\big) & \Longrightarrow & \Big(u(t,x)\ge p^+(t,x_1)-\displaystyle{\frac{\delta}{2}}\Big).\eaa\right.
\ee\par
For all $\xi\ge 2A$ and $x_1-\xi_t\le A$, one has
\be\label{vxi}
v^{\xi}(t,x):=v(t,x_1-\xi,x')\ge p^+(t,x_1-\xi)-\frac{\delta}{2}\ge p^+(t,x_1)-\delta\ge p^-(t,x_1)+\delta
\ee
because $p^+$ is nonincreasing in $x_1$ and one can assume, without loss of generality, that $0<2\delta\le\kappa$, under the notation used in (\ref{fdecreasing}) and (\ref{kappa}). 

\begin{lem}\label{lem4}
For all $\xi\ge 2A$, there holds
\be\label{claim1lem4}
v^{\xi}(t,x)\ge u(t,x)\hbox{ for all }(t,x)\in\R\times\R^N\hbox{ such that }x_1-\xi_t\ge A
\ee
and
\be\label{claim2lem4}
v^{\xi}(t,x)\ge u(t,x)\hbox{ for all }(t,x)\in\R\times\R^N\hbox{ such that }x_1-\xi_t\le A.
\ee
\end{lem}

\noindent{\bf{Proof.}} Fix any $\xi\ge 2A$. We will just prove property (\ref{claim1lem4}), the proof of the second one being similar. Since $u$ is bounded, the nonnegative real number
$$\epsilon^*=\inf\big\{\epsilon>0;\ v^{\xi}(t,x)\ge u(t,x)-\epsilon\hbox{ for all }(t,x)\in\R\times\R^N\hbox{ with }x_1-\xi_t\ge A\big\}$$
is well-defined. Observe that
\be\label{vxieps*}
v^{\xi}(t,x)\ge u(t,x)-\epsilon^*\hbox{ for all }(t,x)\in\R\times\R^N\hbox{ with }x_1-\xi_t\ge A.
\ee
Assume by contradiction that $\epsilon^*>0$. Then there exist a sequence $(\epsilon_n)_{n\in\N}$ of positive real numbers and a sequence $(t_n,x_n)_{n\in\N}=(t_n,x_{1,n},x'_n)_{n\in\N}$ in $\R\times\R^N$ such that
\be\label{epsn}
\epsilon_n\to\epsilon^*\hbox{ as }n\to+\infty,\hbox{ and }x_{1,n}-\xi_{t_n}\ge A,\ v^{\xi}(t_n,x_n)<u(t_n,x_n)-\epsilon_n\hbox{ for all }n\in\N.
\ee
Since $v^{\xi}(t,x)\ge u(t,x)$ when $x_1-\xi_t=A$ from (\ref{Ater}) and (\ref{vxi}), and since $u$ is globally $C^1(\R\times\R^N)$, there exists $\kappa>0$ such that
\be\label{vxieps}
v^{\xi}(t,x)\ge u(t,x)-\frac{\epsilon^*}{2}\hbox{ for all }(t,x)\in\R\times\R^N\hbox{ such that }|x_1-\xi_t-A|<\kappa.
\ee
In particular, there holds
$$x_{1,n}-\xi_{t_n}\ge A+\kappa\hbox{ for large }n.$$
Furthermore, we claim that the sequence $(x_{1,n}-\xi_{t_n})_{n\in\N}$ is bounded. Otherwise, up to extraction of a subsequence, it would converge to $+\infty$. Thus,
$$v^{\xi}(t_n,x_n)-p^-(t_n,x_{1,n}-\xi)=u(t_n,x_{1,n}-\xi,x'_n+\theta)-p^-(t_n,x_{1,n}-\xi)\to 0\hbox{ as }n\to+\infty$$
and
$$u(t_n,x_n)-p^-(t_n,x_{1,n})\to 0\hbox{ as }n\to+\infty.$$
Since $\xi\ge 0$ and $p^-(t,x_1)$ is nonincreasing with respect to $x_1$, it would then follow that
$$\liminf_{n\to+\infty}v^{\xi}(t_n,x_n)-u(t_n,x_n)\ge 0,$$
which contradicts (\ref{epsn}). Thus, the sequence $(x_{1,n}-\xi_{t_n})_{n\in\N}$ is bounded.\par
Remember now that, because of (\ref{xisigma}) and Definition~\ref{def1}, the function $t\mapsto\xi_t$ can be assumed to be uniformly continuous. In particular, the sequence $(\xi_{t_n}-\xi_{t_n-1})_{n\in\N}$ is bounded, whence the sequence $(x_{1,n}-\xi_{t_n-1})_{n\in\N}$ is bounded as well. Moreover, there exists a real number $\rho$ such that
\be\label{defrho3}
0<\rho\le\frac{\kappa}{4}\hbox{ and }|\xi_s-\xi_{s'}|\le\frac{\kappa}{2}\hbox{ for all }(s,s')\in\R^2\hbox{ such that }|s-s'|\le\rho.
\ee
Choose now $K\in\N\backslash\{0\}$ such that
\be\label{defK2}
K\,\rho\ge\max\Big(1,\sup\big\{|x_{1,n}-\xi_{t_n-1}-A|;\ n\in\N\big\}\Big).
\ee
For each $n\in\N$ and $i=0,\ldots,K$, set
$$\tilde{x}_{n,i}=x_{1,n}+\frac{i}{K}\,(\xi_{t_n-1}+A-x_{1,n})$$
and
$$E_{n,i}=\Big[t_n-\frac{i+1}{K},t_n-\frac{i}{K}\Big]\times[\tilde{x}_{n,i}-2\,\rho,\tilde{x}_{n,i}+2\,\rho]\times\big\{x'\in\R^{N-1};\ |x'-x'_n]\le 1\}.$$
Observe that $|\tilde{x}_{n,i+1}-\tilde{x}_{n,i}|\le\rho$ for all $0\le i\le K-1$, from (\ref{defK2}). Furthermore, since $x_{1,n}-\xi_{t_n}>A+\kappa$ for large $n$, say for $n\ge n_0$, it follows from (\ref{defrho3}) and (\ref{defK2}) that
$$\tilde{x}_{n,0}-2\,\rho=x_{1,n}-2\,\rho\ge\xi_t+A\ \hbox{ for all }t_n-\frac{1}{K}\le t\le t_n$$
and for all $n\ge n_0$. Consequently,
\be\label{En0}
E_{n,0}\subset\big\{(t,x)\in\R\times\R^N;\ x_1-\xi_t\ge A\}\ \hbox{ for all }n\ge n_0.
\ee
Thus,
$$w:=v^{\xi}-(u-\epsilon^*)\ge 0\hbox{ in }E_{n,0}$$
and
$$u(t,x)-\epsilon^*<u(t,x)\le p^-(t,x_1)+\delta\hbox{ in }E_{n,0}$$
for all $n\ge n_0$ from (\ref{Ater}) and (\ref{vxieps*}). Since $f(t,x_1,\cdot)$ is nonincreasing in $(-\infty,p^-(t,x_1)+\delta]$, it follows that $u-\epsilon^*$ is a subsolution of (\ref{eq}) in $E_{n,0}$ for all $n\ge n_0$, while $v^{\xi}$ is a supersolution of~(\ref{eq}) in $\R\times\R^N$, because $A$ and $q$ only depend on $t$, and $f(t,x_1,s)$ is nonincreasing in $x_1$.\par
Finally, for all $n\ge n_0$, the globally $C^1(\R\times\R^N)$ function $w$ is nonnegative in~$E_{n,0}$, it satisfies inequations of the type
$$w_t\ge\nabla_x\cdot(A(t)\nabla_xw)+q(t)\cdot\nabla_xw+b(t,x)w\hbox{ in }E_{n,0}$$
where the sequence $(\|b\|_{L^{\infty}(E_{n,0})})_{n\in\N}$ is bounded. Since $w(t_n,\tilde{x}_{n,0},x'_n)=w(t_n,x_n)\to 0$ as $n\to+\infty$, one finally concludes from the linear parabolic estimates that
\be\label{w1n}
w\Big(t_n-\frac{1}{K},\tilde{x}_{n,1},x'_n\Big)\to 0\hbox{ as }n\to+\infty.
\ee
But since
$$\tilde{x}_{n,1}-\xi_{t_n-1/K}\ge\tilde{x}_{n,0}-\rho-\xi_{t_n-1/K}\ge A$$
from (\ref{En0}) for all $n\ge n_0$, it follows from (\ref{vxieps}) and (\ref{w1n}) that $\tilde{x}_{n,1}-\xi_{t_n-1/K}\ge A+\kappa$ for $n$ large enough. By repeating the arguments inductively, one concludes that
$$\tilde{x}_{n,i}-\xi_{t_n-i/K}\ge A+\kappa\hbox{ for all }i=1,\ldots,K\hbox{ and for }n\hbox{ large enough}.$$
One gets a contradiction at $i=K$, since $\tilde{x}_{n,K}=\xi_{t_n-1}+A$.\par
As a conclusion, the assumption $\epsilon^*>0$ was false. Hence, the claim (\ref{claim1lem4}) is proved, and, as already emphasized, the proof of (\ref{claim2lem4}) follows the same scheme.\hfill$\Box$\break

\noindent{\bf{End of the proof of Theorem~\ref{th3}.}} Lemma~\ref{lem4} yields
$$v^{\xi}\ge u\hbox{ in }\R\times\R^N\hbox{ for all }\xi\ge 2A.$$
Now define
$$\xi^*\ =\ \inf\big\{\xi>0,\ v^{\xi'}\ge u\hbox{ in }\R\times\R^N\hbox{ for all }\xi'\ge\xi\big\}.$$
One has $0\le\xi^*\le 2A$, and $v^{\xi^*}(t,x)\ge u(t,x)$ for all $(t,x)\in\R\times\R^N$. Assume now that $\xi^*>0$. Two cases may occur:\par
{\it{Case 1:}} assume here that
$$\inf\big\{v^{\xi^*}(t,x)-u(t,x);\ |x_1-\xi_t|\le A\big\}>0.$$
From the boundedness of $u_{x_1}$, there exists then $\eta_0\in(0,\xi^*)$ such that
\be\label{vxi*}
v^{\xi^*-\eta}(t,x)\ge u(t,x)\hbox{ for all }\eta\in[0,\eta_0]\hbox{ and }|x_1-\xi_t|\le A.
\ee
Since
$$v^{\xi^*}(t,x)\ge u(t,x)\ge p^+(t,x_1)-\frac{\delta}{2}\hbox{ for all }x_1-\xi_t\le-A,$$
one can assume that $\eta_0>0$ is small enough so that
$$v^{\xi^*-\eta}(t,x)\ge p^+(t,x_1)-\delta\hbox{ for all }x_1-\xi_t\le-A.$$
Applying again the arguments used in Lemma \ref{lem4}, one then concludes that, for all $\eta\in[0,\eta_0]$, there holds $v^{\xi^*-\eta}(t,x)\ge u(t,x)$ for all $(t,x)\in\R\times\R^N$ such that $x_1-\xi_t\le -A$ or $x_1-\xi_t\ge A$. Eventually, together with (\ref{vxi*}),
$$v^{\xi^*-\eta}\ge u\hbox{ in }\R\times\R^N$$
for all $\eta\in[0,\eta_0]$, which contradicts the minimality of $\xi^*$. Thus, case 1 is ruled out.\par
{\it{Case 2:}} one then has
$$\inf\big\{v^{\xi^*}(t,x)-u(t,x);\ |x_1-\xi_t|\le A\big\}=0.$$
There exists then a sequence $(t_n,x_n)_{n\in\N}=(t_n,x_{1,n},x'_n)_{n\in\N}$ in $\R\times\R^N$ such that
$$\left\{\baa{l}
|x_{1,n}-\xi_{t_n}|\le A\ \hbox{ for all }n\in\N,\vspace{5pt}\\
u(t_n,x_{1,n}-\xi^*,x'_n+\theta)-u(t_n,x_n)=v^{\xi^*}(t_n,x_n)-u(t_n,x_n)\to 0\ \hbox{ as }n\to+\infty.\eaa\right.$$\par
Fix now any $\sigma>0$ and $m\in\N\backslash\{0\}$. Since $v^{\xi^*}\ge u$ and $v^{\xi^*}$ is a supersolution of~(\ref{eq}) in~$\R\times\R^N$, the linear parabolic estimates then imply that
$$\baa{rcl}
u\left(t_n-\displaystyle{\frac{\sigma}{m}},x_{1,n}-2\xi^*,x'_n+2\theta\right)-u\left(t_n-\displaystyle{\frac{\sigma}{m}},x_{1,n}-\xi^*,x'_n+\theta\right) & & \vspace{5pt}\\
=v^{\xi^*}\left(t_n-\displaystyle{\frac{\sigma}{m}},x_{1,n}-\xi^*,x'_n+\theta\right)-u\left(t_n-\displaystyle{\frac{\sigma}{m}},x_{1,n}-\xi^*,x'_n+\theta\right) & \longrightarrow & 0\ \hbox{ as }n\to+\infty.\eaa$$
By immediate induction, one gets that
$$u\left(t_n-k\frac{\sigma}{m},x_{1,n}-(k+1)\xi^*,x'_n+(k+1)\theta\right)-u\left(t_n-k\frac{\sigma}{m},x_{1,n}-k\xi^*,x'_n+k\theta\right)\mathop{\longrightarrow}_{n\to+\infty}0,$$
for each $k=1,\ldots,m$. Therefore,
$$\limsup_{n\to+\infty}\big|u(t_n-\sigma,x_{1,n}-(m+1)\xi^*,x'_n+(m+1)\theta)-u(t_n,x_n)\big|\le\sigma\|u_t\|_{L^{\infty}(\R\times\R^N)}.$$
Similarly, by considering the points $(t_n-k\sigma/m,x_{1,n}+k\xi^*,x'_n+k\theta)$, one gets that
$$\limsup_{n\to+\infty}\big|u(t_n-\sigma,x_{1,n}+(m-1)\xi^*,x'_n+(m-1)\theta)-u(t_n,x_n)\big|\le\sigma\|u_t\|_{L^{\infty}(\R\times\R^N)}.$$
Hence,
\be\label{limsupbis}\baa{lcl}
\displaystyle{\mathop{\limsup}_{n\to+\infty}}\ \big|u(t_n-\sigma,x_{1,n}-(m+1)\xi^*,x'_n+(m+1)\theta)\ & & \\
\qquad\qquad\qquad-u(t_n-\sigma,x_{1,n}+(m-1)\xi^*,x'_n+(m-1)\theta)\big| & \le & 2\sigma\|u_t\|_{L^{\infty}(\R\times\R^N)}.\eaa
\ee\par
Choose now $\sigma>0$ such that
\be\label{sigmabis}
2\sigma\|u_t\|_{L^{\infty}(\R\times\R^N)}\le\frac{\kappa}{4}.
\ee
But $|x_{1,n}-\xi_{t_n}|\le A$ for all $n\in\N$ and the sequence $(\xi_{t_n}-\xi_{t_n-\sigma})_{n\in\N}$ is bounded from the assumption made in Theorem~\ref{th3}. Therefore, the sequence $(x_{1,n}-\xi_{t_n-\sigma})_{n\in\N}$ is bounded. Let $C\ge 0$ be such that
$$\left\{\baa{lcl}
\big(x_1\ge\xi_t+C\big) & \Longrightarrow & \Big(u(t,x)\le p^-(t,x_1)+\displaystyle{\frac{\kappa}{4}}\Big)\vspace{5pt}\\
\big(x_1\le\xi_t-C\big) & \Longrightarrow & \Big(u(t,x)\ge p^+(t,x_1)-\displaystyle{\frac{\kappa}{4}}\Big).\eaa\right.$$
Since $\xi^*$ is assumed to be positive, there exists $m\in\N\backslash\{0\}$ such that
$$x_{1,n}+(m-1)\xi^*\ge\xi_{t_n-\sigma}+C\hbox{ and }x_{1,n}-(m+1)\xi^*\le\xi_{t_n-\sigma}-C\hbox{ for all }n\in\N.$$
Thus,
$$\baa{rcl}
u(t_n-\sigma,x_{1,n}+(m-1)\xi^*,x'_n+(m-1)\theta) & \le & p^-(t_n-\sigma,x_{1,n}+(m-1)\xi^*)+\displaystyle{\frac{\kappa}{4}}\vspace{5pt}\\
& \le & p^-(t_n-\sigma,x_{1,n})+\displaystyle{\frac{\kappa}{4}}\eaa$$
and
$$\baa{rcl}
u(t_n-\sigma,x_{1,n}-(m+1)\xi^*,x'_n+(m+1)\theta) & \ge & p^+(t_n-\sigma,x_{1,n}-(m+1)\xi^*)-\displaystyle{\frac{\kappa}{4}}\vspace{5pt}\\
& \ge & p^+(t_n-\sigma,x_{1,n})-\displaystyle{\frac{\kappa}{4}}\eaa$$
for all $n\in\N$, because $p^{\pm}$ are nonincreasing in $x_1$. Hence,
$$\baa{c}
u(t_n-\sigma,x_{1,n}-(m+1)\xi^*,x'_n+(m+1)\theta)-u(t_n-\sigma,x_{1,n}+(m-1)\xi^*,x'_n+(m-1)\theta)\vspace{5pt}\\
\ge p^+(t_n-\sigma,x_{1,n})-p^-(t_n-\sigma,x_{1,n})-\displaystyle{\frac{\kappa}{2}}\ge\displaystyle{\frac{\kappa}{2}}\ \hbox{ for all }n\in\N,\eaa$$
by definition of $\kappa$. Therefore,
$$\baa{l}
\displaystyle{\mathop{\limsup}_{n\to+\infty}}\big|u(t_n-\sigma,x_{1,n}-(m+1)\xi^*,x'_n+(m+1)\theta)\vspace{5pt}\\
\qquad\qquad\qquad-u(t_n-\sigma,x_{1,n}+(m-1)\xi^*,x'_n+(m-1)\theta)\big|\ge\displaystyle{\frac{\kappa}{2}},\eaa$$
while it is less than or equal to $\kappa/4$ from (\ref{limsupbis}) and (\ref{sigmabis}).\par
One has then reached a contradiction, which means that $\xi^*=0$. Then,
$$v(t,x_1-\xi,x'+\theta)\ge u(t,x_1,x')\hbox{ for all }(t,x_1,x')\in\R\times\R^N,\ \xi\ge 0\hbox{ and }\theta\in\R^{N-1}.$$
As a consequence, $u$ is nonincreasing in $x_1$ and it does not depend on $x'$. Furthermore, the strong parabolic maximum principle, together with the same arguments as above, implies that $u$ is actually decreasing in $x_1$. That completes the proof of Theorem \ref{th3}.\hfill$\Box$\break

\noindent{\bf{Proof of Theorem~\ref{statio}.}} Assume that all assumptions made in Theorem~\ref{statio} are satisfied. Up to rotation of the frame, one can assume without loss of generality that $e=e_1=(1,0,\ldots,0)$. Consider first the case where $c>0$. There exists $\epsilon\in\{-1,1\}$ such that
$$\frac{\xi_t-\xi_s}{t-s}\to\epsilon\,c\hbox{ as }t-s\to\pm\infty\hbox{ and }\sup\big\{|\xi_t-\epsilon\,c\,t|;\ t\in\R\big\}<+\infty.$$
The function
$$v(t,x)=u(t,x+\epsilon\,c\,t\,e)=u(t,x_1+\epsilon\,c\,t,x')$$
is well-defined for all $(t,x)\in\R\times\overline{\Omega}$ (because $\Omega$ is invariant in the direction $e$) and it satisfies
$$\left\{\baa{rcl}
v_t & = & \nabla_x\cdot(A(x')\nabla_xv)+q(x')\cdot\nabla_xv+\epsilon\, c\,v_{x_1}+f(x',v)\hbox{ in }\R\times\Omega,\vspace{5pt}\\
\mu(x')\cdot\nabla_xv & = & 0\hbox{ on }\R\times\partial\Omega,\eaa\right.$$
because $A$, $q$, $\mu$ and $f$ are independent of $x_1$ (and of $t$). Furthermore, since $p^{\pm}$ only depend on $x'$, $v$ is a transition front connecting $p^-$ and $p^+$, with the sets
$$\tilde{\Omega}^{\pm}_t=\big\{x\in\Omega,\ \pm x_1<0\big\}\hbox{ and }\tilde{\Gamma}_t=\big\{x\in\Omega,\ x_1=0\big\}.$$
With the same type of arguments as in the proof of Theorem \ref{th3} above, one can then fix any $\zeta\in\R$ and slide $v(t+\zeta,x)$ with respect to $v$ in the $x_1$-direction. It follows then that
$$v(t+\zeta,x_1-\xi,x')\ge v(t,x_1,x')\hbox{ for all }(t,x)\in\R\times\overline{\Omega},\ \xi\ge 0\hbox{ and }\zeta\in\R.$$
Therefore, $v$ is independent of $t$ and it is nonincreasing in $x_1$. As above, $v$ is then decreasing in $x_1$. That gives the required conclusion in the case where $c>0$.\par
In the case where $c=0$, the function $t\mapsto\xi_t$ is then bounded. Because of Definition~\ref{def1}, one can then assume, without loss of generality, that $\xi_t=0$ for all $t\in\R$. The functions~$p^{\pm}$ and $f$ may depend on $x_1$, but are assumed to be nonincreasing in $x_1$. For any $\zeta\in\R$ and $\xi\ge 0$, the function $u(t+\zeta,x_1-\xi,x')$ is then a supersolution of the equation (\ref{eq}) which is satisfied by $u$. One can then slide $u(t+\zeta,x_1,x')$ with respect to $u$ in the (positive) $x_1$-direction, and it follows as in the proof of Theorem \ref{th3} that
$$u(t+\zeta,x_1-\xi,x')\ge u(t,x_1,x')\hbox{ for all }(t,x)\in\R\times\overline{\Omega},\ \xi\ge 0\hbox{ and }\zeta\in\R.$$
As usual, one concludes that $u$ does not depend on $t$ and is decreasing in $x_1$.\hfill$\Box$


\end{document}